\def\a             {\alpha}
\def\Ad            {{\mathrm{Ad}}}

\def\Aut           {{\mathrm{Aut}}}
\def\be            {\begin{equation}}
\def\bbC           {\mathbb{C}}

\def\bbR           {\mathbb{R}}
\def\bbT           {\mathbb{T}}
\def\bbZ           {\mathbb{Z}}
\def\bfe           {\mathbf{1}}
\def\can           {\gamma}
\def\canr          {\theta}

\def\cA            {{\mathcal{A}}}
\def\cB            {{\mathcal{B}}}

\def\cD            {{\mathcal{D}}}
\def\cE            {{\mathcal{E}}}

\def\cH            {{\mathcal{H}}}

\def\cM            {{\mathcal{M}}}
\def\cN            {{\mathcal{N}}}
\def\cO            {{\mathcal{O}}}

\def\cS            {{\mathcal{S}}}
\def\cT            {{\mathcal{T}}}
\def\cV            {{\mathcal{V}}}

\def\cZ            {{\mathcal{Z}}}
\newcommand\co[1]  {\overline{{#1}}}

\def\E             {{\mathrm{e}}}
\def\ee            {\end{equation}}
\def\End           {{\mathrm{End}}}
\def\eps           {\varepsilon}

\newcommand\erf[1] {Eq.\ \nolinebreak (\ref{#1})}

\def\ext           {{\mathrm{ext}}}
\newcommand\fig[1] {Fig.\ \nolinebreak \ref{#1}}
\def\Gtwo          {{\mathrm{G}}_2}
\def\Hom           {{\mathrm{Hom}}}
\def\I             {{\mathrm{i}}}
\def\id            {{\mathrm{id}}}
\def\Ind           {{\mathrm{Ind}}}

\def\la            {\lambda}
\def\lan           {\langle}

\def\LSOf          {{\mathit{LSO}}(5)}
\def\LSUn          {{\mathit{LSU}}(n)}
\def\LSUz          {{\mathit{LSU}}(2)}
\def\LISUn         {{\mathit{L}}_I{\mathit{SU}}(n)}
\def\LISUz         {{\mathit{L}}_I{\mathit{SU}}(2)}
\def\LISOf         {{\mathit{L}}_I{\mathit{SO}}(5)}
\def\LIcSUn        {{\mathit{L}}_{I^\rmc}{\mathit{SU}}(n)}
\def\Mat           {{\mathrm{Mat}}}
\def\Mor           {{\mathrm{Mor}}}

\def\MXN           {{}_M {\mathcal{X}}_N}
\def\MXM           {{}_M {\mathcal{X}}_M}
\def\MXMa          {{}_M^{} {\mathcal{X}}_M^\a}
\def\MXMo          {{}_M^{} {\mathcal{X}}_M^0}
\def\MXMop         {{}_{M_+}^{}\!\! {\mathcal{X}}_{M_+}^0}
\def\MXMom         {{}_{M_-}^{}\!\! {\mathcal{X}}_{M_-}^0}
\def\MXMopm        {{}_{M_\pm}^{}\!\! {\mathcal{X}}_{M_\pm}^0}
\def\MXMp          {{}_M^{} {\mathcal{X}}_M^+}
\def\MXMm          {{}_M^{} {\mathcal{X}}_M^-}
\def\MXMpm         {{}_M^{} {\mathcal{X}}_M^\pm}
\def\MXMpp         {{}_{M_+}^{} {\mathcal{X}}_{M_+}^+}
\def\MXMmm         {{}_{M_-}^{} {\mathcal{X}}_{M_-}^-}
\def\MXMppm        {{}_{M_+}^{} {\mathcal{X}}_{M_+}^\pm}
\def\MXMpmpm       {{}_{M_\pm}^{} {\mathcal{X}}_{M_\pm}^\pm}
\def\NXN           {{}_N {\mathcal{X}}_N}
\def\NXNd          {{}_N^{} {\mathcal{X}}_N^{\mathrm{deg}}}

\def\om            {\omega}

\def\op            {{\mathrm{opp}}}

\def\ran           {\rangle}
\def\rmA           {{\mathrm{A}}}
\def\rmD           {{\mathrm{D}}}
\def\rmE           {{\mathrm{E}}}
\def\rmb           {{\mathrm{b}}}
\def\rmv           {{\mathrm{v}}}
\def\rms           {{\mathrm{s}}}
\def\rmc           {{\mathrm{c}}}
\def\sig           {\sigma}

\def\SLZ           {{\mathit{SL}}(2;\bbZ)}
\def\SOf           {{\mathit{SO}}(5)}
\def\SOn           {{\mathit{SO}}(n)}

\def\SUd           {{\mathit{SU}}(3)}

\def\SUn           {{\mathit{SU}}(n)}

\def\SUz           {{\mathit{SU}}(2)}

\def\SUf           {{\mathit{SU}}(4)}
\def\sys           {\Delta}
\newcommand\tmat[1]{{}^{\mathrm{t}} {#1}}
\def\tn            {{\tilde{n}}}
\def\tr            {{\mathrm{tr}}}

\def\typei         {type \nolinebreak I}
\def\typeii        {type \nolinebreak II}
\def\typeiii       {type \nolinebreak III}

\def\thinlines{\allinethickness{0.3pt}}
\def\thicklines{\allinethickness{1.0pt}}
\def\Thicklines{\allinethickness{2.0pt}}

\documentclass{conm-p-l}
\usepackage{amssymb,amsfonts,latexsym,epic,eepic}

\theoremstyle{definition}

\theoremstyle{remark}

\numberwithin{equation}{section}

\begin{document}

\title{Modular Invariants from Subfactors}

\author{Jens B\"ockenhauer}
\address{School of Mathematics, University of Wales Cardiff,
PO Box 926, Senghennydd Road, Cardiff CF24 4YH, Wales, U.K.}
\email{BockenhauerJM@cardiff.ac.uk}
\thanks{This project was supported by the EU TMR Network in
Non-Commutative Geometry.}

\author{David E.\ Evans}
\address{School of Mathematics, University of Wales Cardiff,
PO Box 926, Senghennydd Road, Cardiff CF24 4YH, Wales, U.K.}
\email{EvansDE@cardiff.ac.uk}
\thanks{Lectures given by the second author at
``Quantum Symmetries in Theoretical Physics and Mathematics'',
10--21 January 2000, Bariloche, Argentine.}


\date{June 1, 2000 and, in revised form, June 22, 2024.}


\begin{abstract}
In these lectures we explain the intimate relationship
between modular invariants in conformal field theory
and braided subfactors in operator algebras.
A subfactor with a braiding determines a matrix $Z$
which is obtained as a coupling matrix comparing two
kinds of braided sector induction (``$\a$-induction'').
It has non-negative integer entries, is normalized
and commutes with the S- and T-matrices arising from
the braiding. Thus it is a physical modular invariant
in the usual sense of rational conformal field theory.
The algebraic treatment of conformal field theory models,
e.g.\ $\SUn_k$ models, produces subfactors which realize
their known modular invariants.
Several properties of modular invariants have so far
been noticed empirically and considered mysterious
such as their intimate relationship to graphs,
as for example the A-D-E classification for $\SUz_k$.
In the subfactor context these properties can be
rigorously derived in a very general setting.
Moreover the fusion rule isomorphism for maximally
extended chiral algebras due to Moore-Seiberg,
Dijkgraaf-Verlinde finds a clear and very general
proof and interpretation through intermediate subfactors,
not even referring to modularity of $S$ and $T$.
Finally we give an overview on the current state of
affairs concerning the relations between the
classifications of braided subfactors and two-dimensional
conformal field theories. We demonstrate in particular
how to realize twisted (\typeii) descendant modular invariants
of conformal inclusions from subfactors and illustrate the
method by new examples.
\end{abstract}

\maketitle


\section{Introduction and overview}

A subfactor in its simplest guise arises from a group
action $M^G\subset M$, the fixed point algebra $M^G$
in the ambient von Neumann algebra $M$ where a group
$G$ acts upon. If say the group is finite and acts
outerly on $M$ (equivalently $(M^G)'\cap M=\bbC \bfe$,
where the prime denotes the commutant) and both the
group and the algebra are amenable, then we can recover
both the group and the action from the inclusion
$M^G\subset M$. (If $M$ is not amenable, i.e.\ hyperfinite,
one may recover the group but not the action as in free
group factors in free probability theory).
However we will concentrate on (infinite-dimensional)
hyperfinite von Neumann algebras $M$
which are inductive limits of finite dimensional
algebras and are factors i.e.\ have trivial center
$M'\cap M=\bbC\bfe$. A subfactor $N\subset M$ is then
an inclusion of one factor in another, which is thought
to represent a deformation of a group, for us we will
restrict to the case where we only think of those
inclusions which are deviants of finite groups.
(Cf.\ \cite{EK} as a general reference.)

\begin{figure}[htb]
\begin{center}
\unitlength 0.15mm
\begin{picture}(600,230)
\thinlines
\multiput(0,0)(40,0){6}{\circle*{5}}
\multiput(20,40)(40,0){5}{\circle*{5}}
\thicklines
\multiput(0,0)(40,0){5}{\line(1,2){20}}
\multiput(20,40)(40,0){5}{\line(1,-2){20}}
\thinlines
\multiput(300,0)(60,0){6}{\circle*{5}}
\multiput(330,40)(60,0){5}{\circle*{5}}
\multiput(360,80)(60,0){4}{\circle*{5}}
\multiput(390,120)(60,0){3}{\circle*{5}}
\multiput(420,160)(60,0){2}{\circle*{5}}
\put(450,200){\circle*{5}}
\thicklines
\path(300,0)(450,200)(600,0)(540,0)(420,160)(480,160)
(360,0)(330,40)(570,40)(540,0)(420,0)(360,80)(540,80)
(480,0)(390,120)(510,120)(420,0)(300,0)
\Thicklines
\multiput(330,0)(60,0){5}{\vector(1,0){0}}
\multiput(360,40)(60,0){4}{\vector(1,0){0}}
\multiput(390,80)(60,0){3}{\vector(1,0){0}}
\multiput(420,120)(60,0){2}{\vector(1,0){0}}
\multiput(450,160)(60,0){1}{\vector(1,0){0}}
\multiput(315,20)(30,40){5}{\vector(-3,-4){0}}
\multiput(375,20)(30,40){4}{\vector(-3,-4){0}}
\multiput(435,20)(30,40){3}{\vector(-3,-4){0}}
\multiput(495,20)(30,40){2}{\vector(-3,-4){0}}
\multiput(555,20)(30,40){1}{\vector(-3,-4){0}}
\multiput(585,20)(-30,40){5}{\vector(-3,4){0}}
\multiput(525,20)(-30,40){4}{\vector(-3,4){0}}
\multiput(465,20)(-30,40){3}{\vector(-3,4){0}}
\multiput(405,20)(-30,40){2}{\vector(-3,4){0}}
\multiput(345,20)(-30,40){1}{\vector(-3,4){0}}
\end{picture}
\end{center}
\caption{Fusion graphs of fundamental generators $\Box$
         of systems for $\SUz_{10}$ and $\SUd_5$}
\label{agraphs}
\end{figure}
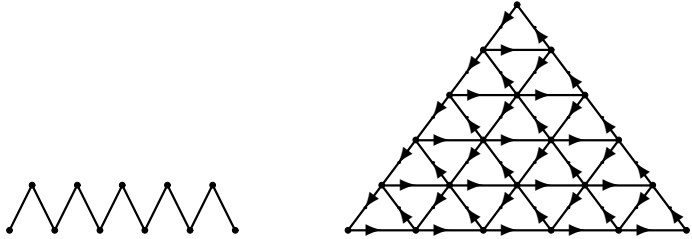

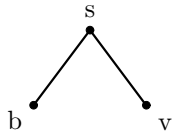
\begin{figure}[htb]
\begin{center}
\unitlength 0.5mm
\begin{picture}(40,40)
\thinlines
\multiput(5,5)(30,0){2}{\circle*{2}}
\put(20,25){\circle*{2}}
\thicklines
\path(5,5)(20,25)(35,5)
\put(0,1){\makebox(0,0){$\rmb$}}
\put(20,30){\makebox(0,0){$\rms$}}
\put(40,0){\makebox(0,0){$\rmv$}}
\end{picture}
\end{center}
\caption{Dynkin diagram A$_3$ as fusion graph}
\label{A3}
\end{figure}

Rather than a group of $\ast$-automorphisms of
a von Neumann algebra $M$, we will more generally consider
a system $\sys$ of $\ast$-endomorphisms which is closed
under composition
\[
\la \circ \mu = \bigoplus_{\nu\in\sys} N_{\la,\mu}^\nu\,\nu
\]
for a suitable notion of addition of endomorphisms (for which we
will need infinite von Neumann factors and consider endomorphisms
up to inner equivalence, i.e.\ as sectors \cite{L2}) and
non-negative integral coefficients $N_{\la,\mu}^\nu$. In our
relationship with modular invariant partition functions in
conformal field theory, our starting point will be a system of
endomorphisms labelled by vertices of graphs as e.g.\ given in
\fig{agraphs}.
Each $\la\in\sys$ defines a matrix
$N_\la=[N_{\la,\mu}^\nu]_{\mu,\nu}$
of multiplication by $\la$ so that in the above setting the
graph of $N_\Box$ where $\Box$ is the fundamental generator
is as described in the figures. For example in the case of
the Dynkin diagram A$_3$, as in \fig{A3}.
Here we labelled the vertices by b, s, v, and the graph
represents the `fusion' by s, and so the multiplication
by s gives the sum of nearest neighbors:
\[ \rms \cdot \rmb = \rms \,,\qquad
\rms \cdot \rms = \rmb \oplus \rmv \,,\qquad
\rms \cdot \rmv = \rms \,. \]
(Here and in general it is understood that an unoriented
edge represents an arrow in both directions.)

These are the well-known fusion rules of the
conformal Ising model. A treatment of the Ising model in the
framework of local quantum physics realizing these fusion
rules in terms of endomorphisms on von Neumann factors was
carried out in \cite{B1}, building on \cite{MaSch}.
The transfer matrix formalism allows one to study classical
statistical mechanical models via non-commutative
operator algebras. A study of the Ising model in this
framework was carried out in \cite{AE,EL,CE}.
The fundamental example of this non-commutative framework
for understanding the Ising model was the driving force
towards the present work on understanding modular invariant
partition functions via non-commutative operator algebras
(cf.\ the lecture by the second author at the CBMS meeting
in Eugene, Oregon, September 1993).

Using associativity of the fusion product one obtains
for the fusion matrices
\[ N_\la N_\mu = \sum\nolimits_\nu N_{\la,\mu}^\nu N_\nu \,, \]
i.e.\ the matrices $N_\la$ themselves give a
(``regular'') representation of the fusion rules of $\sys$.
Usually the system $\sys$ will be closed under a certain
conjugation $\la\mapsto\co\la$ (generalizing the notion of
inverse and conjugate representation in a group and group
dual, respectively) which is anti-multiplicative and additive.
This will mean that the transpose of $N_\la$ is $N_{\co\la}$.
If we start with a system obeying commutative fusion rules
(which will not always be the case), the collection
$\{ N_\la \}_{\la\in\sys}$ will therefore constitute a
family of normal commuting matrices, and hence be
simultaneously diagonalizable, with spectra
spec$(N_\la)=\{\gamma^\la_\rho\}_\rho$.
In fact their spectra will be labelled naturally
by the entire set $\sys$ itself, i.e.\ we will
have $\rho\in\sys$. In this diagonalization we have
\begin{equation}\label{1drep}
\gamma^\la_\rho \gamma^\mu_\rho =
\sum\nolimits_\nu N_{\la,\mu}^\nu \gamma^\nu_\rho \,,
\end{equation}
i.e.\ the eigenvalues provide one-dimensional
representations of the fusion rules.
The matrix $\gamma^\la_\rho$ is invertible and we can
invert \erf{1drep} to obtain the Verlinde formula \cite{Ve}
\begin{equation}\label{verlinde}
N_{\la,\mu}^\nu = \sum\nolimits_\rho
\frac{S_{\la,\rho}}{S_{0,\rho}}
S_{\mu,\rho} S_{\nu,\rho}^* \,.
\end{equation}
Here we write the eigenvalues of $N_\la$ as
$\gamma^\la_\rho = S_{\la,\rho}/S_{0,\rho}$,
where the label ``$0$'' refers to the distinguished
identity element (``vacuum'') of the fusion rules,
and $S_{0,\rho}=(\sum_\la |\gamma^\la_\rho|^2)^{1/2}$.
(See \cite{F1} for fusion rules in the context of
conformal field theory.)

In our subfactor approach to modular invariants we will
have representations of the Verlinde fusion rules
appearing naturally, with spectrum a proper subset of $\sys$
and with multiplicities $Z_{\la,\la}$, $\la\in\sys$, given
by the diagonal part of a modular invariant. The
representation matrices can be interpreted a adjacency
matrices of graphs associated with modular invariants.

Modular invariant partition functions arise as
continuum limits in statistical mechanics and play
a fundamental role in conformal field theory.
Recall that a modular invariant partition function
is of the form (cf.\ Zuber's lectures, or see
\cite{DMS,F,K,DiF,G4} for more details on these matters)
\[ Z(\tau) = \sum\nolimits_{\la,\mu} Z_{\la,\mu}
\chi_\la(\tau) \chi_\mu(\tau)^* \,.\]
Here $\chi_\la=\tr(q^{L_0-c/24})$,
$q=\E^{2\pi\I\tau}$, is the trace in the
irreducible representation of a chiral algebra,
which for us will be a positive energy representation
of a loop group with the conformal Hamiltonian $L_0$
being the infinitesimal generator of the rotation group
on the circle. (More typically we would take
un-specialized characters in order to have linearly
independent characters.
See for example \cite{DP} or \cite[Sect.\ 8.3]{EK}
for explicit computations with corner transfer matrices
and derivations of the Virasoro characters in the context
of the Ising model.)
Then the action of the modular group $\SLZ$ on
$q=\E^{2\pi\I\tau}$ via
$\cS=\left({0\atop 1}{-1\atop 0}\right)$:
$\tau\mapsto -1/\tau$, and
$\cT=\left({1\atop 0}{1\atop 1}\right)$:
$\tau\mapsto \tau+1$, transforms the family
of characters $\{\chi_\la\}$ linearly.
More precisely, there are matrices $S$ and $T$
such that
\[
\chi_\la(-1/\tau)= \sum\nolimits_\mu S_{\la,\mu} \chi_\mu(\tau) \,,
\qquad \chi_\la(\tau+1) = \sum\nolimits_\mu T_{\la,\mu} \chi_\mu(\tau) \,.
\]
Note first that what is remarkable about the Verlinde formula,
\erf{verlinde}, is that the matrix which diagonalizes the fusion rules
is the same as the modular matrix $S$ which transforms the characters
(e.g.\ the Kac-Peterson matrix for current algebra models,
see \cite{K,F}).
It is also remarkable that this matrix is symmetric:
$S_{\la,\mu}=S_{\mu,\la}$.

{}From physical considerations we will require solutions to the
matrix equations $ZS=SZ$, $ZT=TZ$, subject to the constraint
$Z_{0,0}=1$ (``uniqueness of the vacuum'') and the
``coupling matrix'' $Z$ having
only non-negative integer entries (from multiplicities
of the representations). There will always be at least one
solution, the diagonal partition function
\[
Z= \sum\nolimits_\la |\chi_\la|^2
\]
(or $Z_{\la,\mu}=\delta_{\la,\mu}$), or more generally
there may be permutation invariants
\[
Z = \sum\nolimits_\la \chi_\la \chi_{\om(\la)}^* \,,
\]
whenever $\om$ is a permutation of the labels which preserves the
fusion rules, the vacuum, and the ``conformal dimensions''. Moore
and Seiberg argue in \cite{MS2} (see also \cite{DV}) that after a
``maximal extension of the chiral algebra'' (the hardest part is
to make this mathematically precise) the partition function of a
RCFT is at most a permutation matrix
$Z^\ext_{\tau,\tau'}=\delta_{\tau,\om(\tau')}$, where $\tau,\tau'$
label the representations of the extended chiral algebra and now
$\om$ denotes a permutation of these with analogous invariance
properties. Decomposing the extended characters $\chi^\ext_\tau$
in terms of the original characters $\chi_\la$, we have
$\chi^\ext_\tau=\sum_\la b_{\tau,\la} \chi_\la$ for some
non-negative integral branching coefficients $b_{\tau,\la}$. The
maximal extension yields the coupling matrix expression
\[ Z_{\la,\mu}=\sum\nolimits_\tau
b_{\tau,\la} b_{\om(\tau),\mu} \,.\] There is a distinction
\cite{DZ2} between so-called \typei\ invariants which arise from
the diagonal invariant of the maximal extension, i.e.\ for which
$\om$ is trivial, and \typeii\ invariants corresponding to
non-trivial automorphisms of the extended fusion rules. The
coupling matrix of a \typei\ invariant is in particular symmetric
whereas \typeii\ invariants need not be so but still the ``vacuum
coupling'' is symmetric: $Z_{0,\la}=Z_{\la,0}$ for all labels
$\la$. To allow more generally for possibly non-symmetric vacuum
coupling, $Z_{0,\la}\neq Z_{\la,0}$, one may need different
extensions for the left and right chiral algebra \cite{BE4}
(see also \cite{KS} where this possibility is explicitly
addressed in the context of simple current extensions), and
then the distinction between \typei\ and \typeii\ modular
invariants does no longer make sense.\footnote{Surprisingly
enough, all known modular invariants of $\SUn_k$ models are
entirely symmetric. Nevertheless there are known modular
invariants of other models with non-symmetric (``heterotic'')
vacuum coupling --- see Section \ref{hetero}.}

A simple argument of Gannon \cite{G1} shows that there are at most
finitely many solutions to our modular invariant problem. Since
$d_\la=S_{\la,0}/S_{0,0}$ will be positive and at least 1 (the
$d_\la$'s will be the Perron-Frobenius weights of the graphs as in
\fig{agraphs}, or indeed $d_\la$ will be the statistical dimension
of $\la$ as a sector in the von Neumann algebra theory or the
square root of the Jones index \cite{J1}), we obtain from
$SZS^*=Z$ that $\sum_{\la,\mu}Z_{\la,\mu}\le\sum_{\la,\mu} d_\la
Z_{\la,\mu} d_\mu = 1/S_{0,0}^2$. Consequently each integer
$Z_{\la,\mu}$ must be bounded by $1/S_{0,0,}^2=\sum_\la d_\la^2$
(from unitarity of the S-matrix), so that there are only finitely
many solutions. Note that this bound will be our ``global index''
$w$, and this suggests a strong relation between Gannon's argument
and Ocneanu's rigidity theorem (presented at a conference in
January 1997 in Madras, India), the latter implying the finiteness
of the number of subequivalent paragroups for a given paragroup.

Gannon's estimate can even be refined to the inequality
\begin{equation}\label{Zbound}
Z_{\la,\mu} \le d_\la d_\mu
\end{equation}
for each individual entry of a modular invariant
coupling matrix as follows. As by Verlinde's
formula, \erf{verlinde}, the eigenvalues of the non-negative
fusion matrices $N_\la$ are given by
$\gamma^\la_\rho=S_{\la,\rho}/S_{0,\rho}$,
Perron-Frobenius theory tells us that
$|\gamma^\la_\rho|$ is bounded by the
Perron-Frobenius eigenvalue $\gamma^\la_0$,
so that
$|S_{\la,\rho}|\le S_{\la,0} S_{0,\rho}/S_{0,0}$
(cf.\ \cite{G4}).
Commutativity of $Z$ with the unitary $S$ then yields
\[
Z_{\la,\mu}=\sum_{\rho,\nu} S_{\la,\rho}Z_{\rho,\nu}S_{\nu,\mu}^*
\le \sum_{\rho,\nu} |S_{\la,\rho}| Z_{\rho,\nu} |S_{\nu,\mu}|
\le d_\la Z_{0,0} d_\mu \,,
\]
which provides \erf{Zbound} by the normalization $Z_{0,0}=1$.

\section{Operator algebraic input}

We will study the classification of modular invariants and
construction of maximal extensions through subfactors,
in particular starting with braided systems of endomorphisms
on loop group factors which are purely infinite factors with
no traces, or more precisely \typeiii$_1$.
Recall that a factor is \typei\ when there is a trace on the
algebra taking discrete values on projections, \typeii\ when
there is a trace that takes continuous values.
(We hope that the classification of factors into types
I, II and III will not be confused with the distinction
of \typei\ and \typeii\ modular invariants --- it has
nothing to do with it.)
A trace on a von Neumann algebra $M$ is a (possibly unbounded)
linear functional $\tau$ satisfying $\tau(ab)=\tau(ba)$,
$a,b\in M$, where the algebra or trace is finite if
$\tau(\bfe)<\infty$, or infinite otherwise. Thus a finite
\typei\ factor is (isomorphic to) $\Mat(n)=\End(\bbC^n)$,
the $n\times n$ complex matrices, the infinite factor is
$B(H)$, the bounded linear operators on an infinite-dimensional
Hilbert space. A factor is of \typeiii\ (or purely infinite)
otherwise, there is no trace and every non-zero projection $p$
is equivalent to the unit in the sense that there is a partial
isometry $v$ in the algebra such that $v^*v=\bfe$ and
$vv^*=p$. The factors relevant for RCFT are amenable, in the
sense that they are hyperfinite, the completions of
unions of finite-dimensional algebras. Murray and von Neumann
showed that there is an unique hyperfinite II$_1$
(i.e.\ finite \typeii) factor which can be realized as
for example the infinite tensor product of matrix algebras
(arbitrarily chosen as long as they are non-commutative)
completed with respect to the trace, i.e.\ use the trace
$\tau$ (constructed as the tensor product of traces over
the matrices) to define an inner product on $\cM$
(the algebraic tensor product)
$\langle a,b \rangle=\tau(b^*a)$.
Letting $\Omega=\bfe$ regarded as a vector in the completion
$\cH$ of $\cM$ with respect to this inner product, we can let
$\cM$ act on $\cH$ by the induced left action of $\cM$ on
itself, and the hyperfinite II$_1$ factor is the
von Neumann algebra generated by $\cM$ in this representation.
There is by Connes \cite{C} an unique hyperfinite II$_\infty$
(i.e.\ infinite \typeii\ factor) which is $R\otimes B(H)$ where
$R$ is the unique hyperfinite \nolinebreak II$_1$ factor and
$B(H)$ is \typei$_\infty$.
There is a finer classification of \typeiii\ factors
into III$_\la$, $0\le \la \le 1$.
For each $\la\in(0,1]$ there is by
Connes an unique hyperfinite \nolinebreak III$_\la$ factor
(the analysis completed by Haagerup \cite{Hp} in the case $\la=1$).
The \typeiii$_0$ factors are classified by their flow of weights.

In the semi-finite case (I or II) where there is a trace
$\tau$, we can define a conjugation $J$, a conjugate
linear map of the Hilbert space $\cH$ of the trace,
by $J:a\mapsto a^*$, $a\in M\subset\cH$, or
$Ja\Omega=a^*\Omega$. Then $J$ is isometric
because $\tau$ is a trace and interchanges left and
right multiplication, indeed $JMJ=M'$.
Thus $M$ and $M'$ are of comparable size.
If $M=\Mat(n)=\bbC^n\otimes\overline{\bbC^n}$
is finite dimensional, then acting on itself
(regarded as a Hilbert space) $M$
becomes $M\otimes\bfe$ with commutant $\bfe\otimes M$.
In general represent a factor $M$ on a Hilbert space $H$
with vector $\Phi\in H$ cyclic for $M$,
i.e.\ $H=\overline{M\Phi}$, which is also cyclic for $M'$,
$H=\overline{M'\Phi}$. (Take a faithful normal state $\varphi$
on $M$ and the associated Hilbert space.)
Then we can define $S:a\Phi\mapsto a^*\Phi$, $a\in M$,
and take the polar decomposition $S=J\Delta^{1/2}$,
where $J$ is a conjugation and $\Delta$ the
(possibly unbounded) Tomita-Takesaki modular operator.
Then Tomita-Takesaki theory \cite{T} tells us that
$JMJ=M'$, and $\sigma_t=\Ad (\Delta^{\I t})$ defines
a one-parameter automorphism group of $M$ which
describes how far the vector state
$\varphi(\cdot)=\langle \cdot \Phi,\Phi\rangle$
is from being a trace,
$\varphi(ab)=\varphi(b\sigma_\I(a))$ for
analytic $a,b\in M$. In the case of a semi-finite algebra,
and if $\varphi$ is a trace, then $S=J$, $\Delta=\bfe$, and
$\sigma_t=\id$, whilst for other choices of cyclic and
separating vectors $\Phi$, the Tomita-Takesaki
modular group $\sigma_t$ is at least inner.

Now consider the case of an infinite subfactor $N\subset M$,
i.e.\ both factors $N$ and $M$ are infinite which means that
they contain isometries with range projections being
different from the identity.
Then we can represent $M$ on a Hilbert space $H$ where
there is a vector $\Phi$ which is cyclic and separating
for both $N$ and $M$. Taking the corresponding
Tomita-Takesaki modular conjugations $J_N$ and $J_M$
where $J_N N J_N=N'$, $J_M M J_M=M'$, we define
\[
\can=\Ad(J_N J_M) |_M : M \rightarrow M' \subset N'
\rightarrow N
\]
called the canonical endomorphism \cite{L1} from $M$ into $N$.
Different choices of Hilbert spaces and cyclic and separating
vectors only amount to a change $\can\rightarrow\Ad(u)\circ\can$
by a unitary $u\in N$, i.e.\ the $N$-$M$ sector determined by
$\can$ is well-defined. (If $\rho\in\Mor(A,B)$ is a unital
morphism from $A$ to $B$, the $B$-$A$ sector $[\rho]$ is the
equivalence class of $\rho$ where $\rho'\simeq\rho$ iff
$\rho'=\Ad(u)\circ\rho$ for unitaries $u\in B$.) We then have an
inclusion of factors:
\begin{equation}\label{Jtt}
\can(N)\subset\can(M)\subset N \subset M
\subset M_1=\Ad(J_MJ_N)(N) \,.
\end{equation}
We can continue upwards (called the Jones tower)
or downwards (called the Jones tunnel) but the sequence
is of period two, e.g.\ the inclusion $\can(N)\subset\can(M)$
is isomorphic to $N\subset M$, and
$\can(M)=\Ad(J_N)(M')\subset N$ is isomorphic
to $M=\Ad(J_M)(M')\subset\Ad(J_MJ_N)(N)=M_1$.
This periodicity reduces to that between a group and its dual
$G\leftrightarrow\hat{G}$ in the case of a group subfactor
tower $M^G\subset M\subset M\rtimes G$.
So there are basically two canonical endomorphisms,
$\can\in\End(M)$ and $\canr\in\End(N)$, where
$\canr=\can|_N$. We call $\can$ the canonical endomorphism,
and $\canr$ the dual canonical endomorphism for $N\subset M$.

The tower can be identified with the Jones extensions,
in the case of finite index obtained by adjoining a
sequence of projections satisfying the Temperley-Lieb
relations. We could define the Jones index using the
Pims\-ner-Popa inequality as follows. If $E:M\rightarrow N$
is a conditional expectation (a projection of norm one),
then let $\Ind(E)$ be the best constant $\xi$ such that
$E(x^*x)\ge \xi^{-1} x^*x$ for all $x\in M$. Then the
(Jones) index $[M:N]$ is the infimum of $\Ind(E)$ over
all expectations $E$, and there is an unique expectation
called the minimal expectation which realizes the index.

For us, all the relative commutants $N'\cap M_j$,
$M'\cap M_j$, in the tower will be finite-dimensional
and moreover
$N'\cap M_j\subset N'\cap M_{j+1}$,
$M'\cap M_j\subset M'\cap M_{j+1}$,
will be described by {\em finite} graphs.
Due to the periodicity of the tower, only two graphs
appear here, the principal and dual principal graph.
The finiteness of the graphs (equivalent to the
finiteness in RCFT) will imply finite index and the
Jones index will be the square of the norm of either graph.

One question which will engage us will be whether a particular
endomorphism of a factor $N$ should be a dual canonical
endomorphism (of some subfactor $N\subset M$ without
any a priori knowledge of what $M$ should be). For example
if $Z$ is a modular invariant, we can consider
$\bigoplus_{\la\in\sys}Z_{0,\la}[\la]$,
$\bigoplus_{\la\in\sys}Z_{\la,0}[\la]$,
$\bigoplus_{\la\in\sys}Z_{\la,\mu}[\la\otimes\mu^\op]$
as candidates for (the sectors of) dual canonical
endomorphisms (on $N$, $N$, $N\otimes N^\op$, respectively,
if $\sys$ is a system of endomorphisms of $N$).

Before we go any further let us formalize the notion of algebraic
operations and sectors. Let $A$ and $B$ be \typeiii\ von Neumann
factors. A unital $\ast$-homomorphism $\rho:A\rightarrow B$ is
called a $B$-$A$ morphism, and we write $\rho\in\Mor(A,B)$. The
positive number $d_\rho=[B:\rho(A)]^{1/2}$ is called the
statistical dimension of $\rho$; here $[B:\rho(A)]$ is the Jones
index of the subfactor $\rho(A)\subset B$. Now if
$\sigma\in\Mor(B,C)$ (with $C$ being another \typeiii\ factor)
then the multiplication or ``fusion''
\[
[\sigma][\rho]=[\sigma\rho]
\]
is well defined on sectors.
(We usually abbreviate $\sigma\rho\equiv\sigma\circ\rho$.)
For $\tau_1,\tau_2\in\Mor(A,B)$ take isometries
$t_1,t_2\in B$ such that $t_1t_1^*+t_2t_2^*=\bfe$
which we can find by infiniteness of B.
Then define the sum
\[
[\tau_1]\oplus[\tau_2]=[\tau] \,, \quad \mbox{where}
\quad \tau(a)=t_1\tau_1(a)t_1^*+t_2\tau_2(a)t_2^* \,,
\quad a\in A \,.
\]
This is well-defined as if $s_1,s_2\in B$ is another
choice of isometries satisfying
$s_1s_1^*+s_2s_2^*=\bfe$ then $u=s_1t_1^*+s_2t_2^*$
is a unitary in $B$, intertwining $\tau$ and $\tau'$
where $\tau'(a)=s_1\tau_1(a)s_1^*+s_2\tau_2(a)s_2^*$
for all $a\in A$. This notion of a sum is basically
writing $[\tau_1]\oplus[\tau_2]$ as a $2\times 2$ matrix
$\left({\tau_1(\cdot)\atop 0}{0\atop \tau_2(\cdot)}\right)$
in $B$ using the infiniteness of $B$ to achieve the
matrix decomposition.
If $\rho$ and $\sig$ are $B$-$A$
morphisms with finite statistical dimensions, then
the vector space of intertwiners
\[ \Hom(\rho,\sig)=\{ t\in B: t\rho(a)=\sig(a)t \,,
\,\, a\in A \}  \]
is finite-dimensional, and we denote its dimension by
$\lan\rho,\sig\ran$.
Note that for $\tau,\tau_1,\tau_2$ and $t_1,t_2$ as
above we have e.g.\ $t_1\in\Hom(\tau_1,\tau)$.
The impossibility of decomposing some $\rho\in\Mor(A,B)$
as $[\rho]=[\rho_1]\oplus[\rho_2]$ for some
$\rho_1,\rho_2\in\Mor(A,B)$,
or irreducibility is then equivalent to the subfactor
$\rho(A)\subset B$ being irreducible, i.e.\
$\rho(A)'\cap B=\bbC\bfe$.

For groups (and group duals) we have a notion of a
conjugate of $\la$, namely the inverse $\la^{-1}$
of $\la$ (respectively the conjugate representation).
There is a similar notion for sectors.
For an irreducible $\la\in\Mor(A,B)$, an irreducible
morphism $\co\la\in\Mor(B,A)$ is a representative
of the conjugate sector if $[\la\co\la]$
or $[\co\la\la]$ contain the identity sector
($[\id_A]$ or $[\id_B]$, respectively), and the
multiplicity is then automatically one for both
cases \cite{I1}. More generally, for an arbitrary
morphism $\rho\in\Mor(A,B)$ of finite statistical
dimension $d_\rho$, an $A$-$B$ morphism $\co\rho$ is a
conjugate morphism if there are isometries
$r_\rho\in\Hom(\id_A,\co\rho\rho)$ and
${\co r}_\rho\in\Hom(\id_B,\rho\co\rho)$ such that
\begin{equation}\label{riso}
\rho(r_\rho)^* {\co r}_\rho=d_\rho^{-1}\bfe_B \qquad\mbox{and}
\qquad \co\rho({\co r}_\rho)^* r_\rho=d_\rho^{-1}\bfe_A \,.
\end{equation}

Recall the tower-tunnel of \erf{Jtt}. Suppose
$E:M\rightarrow N$ is a conditional expectation
and $\varphi$ is a faithful normal state on $N$,
and set $\omega=\varphi\circ E$. Then $\omega$
is a faithful normal state on $M$ such that
$\omega\circ E=\omega$. Take the GNS Hilbert space $H$
of this state on $M$, with cyclic and separating vector
$\Omega$. We can identify this space with
our previous Hilbert space (where there is vector
$\Phi$ being cyclic and separating for both
$N$ and $M$) and the actions coincide. However
$\Phi$ is not identified with $\Omega$, as
$\overline{N\Omega}$ is a proper subspace if
$N\neq M$, with orthogonal
Jones projection
$e_N:\overline{M\Omega}\rightarrow\overline{N\Omega}$
such that $m\Omega\mapsto E(m)\Omega$, $m\in M$.
Define $v':n\Phi\mapsto n\Omega$, $n\in N$, on $H$
so that $v'\in N'$, and $v'{v'}^*=e_N$.
Then $v_1=\Ad(J_M)(v')$ is an isometry in $M_1$,
and also $v_1v_1^*=e_N$. It is easily checked
starting from $J_Mv'=v'J_N$ that $v_1$ is an
intertwiner in $\Hom(\id_{M_1},\can_1)$,
where $\can_1=\Ad(J_{M_1}J_M)$ is the canonical
endomorphism of $M\subset M_1$. Thus, by translating in
the tunnel-tower the canonical and dual canonical
endomorphism contain the identity sector.

Denoting by $\iota:N\hookrightarrow M$ the inclusion
homomorphism we put $\co\iota:M\rightarrow N$ as
$\co\iota(m)=\can(m)$, $m\in M$. Then $\can=\iota\co\iota$
and $\canr=\co\iota\iota$ both contain the identity sector
so that $\co\iota$ is in fact a conjugate morphism for $\iota$.
Similarly, if $\la\in\End(N)$, we can take $\can_\la$,
$\canr_\la$ to be the canonical and dual canonical
endomorphisms of the inclusion $\la(N)\subset N$.
Then we can set $\co\la=\la^{-1}\can_\la$, which is
well defined so that $\la\co\la=\can_\la$.
In the group case, say if we have an outer action
$\alpha : G \rightarrow \Aut(M)$ of a finite group $G$
on a \typeiii\ factor $M$ and let $N$ be the
corresponding fixed point algebra, $N=M^G$,
then $\can$ decomposes as a sector into the group elements,
$[\can]=\bigoplus_{g\in G} [\alpha_g]$, whereas
the decomposition of $\canr$ is according to the
group dual $\hat{G}$, i.e.\
$[\canr]=\bigoplus_{\pi\in\hat{G}} d_\pi [\rho_\pi]$,
with multiplicities given by the dimensions $d_\pi$
of the irreducible representations $\pi$ of $G$.

\begin{figure}[htb]
\begin{center}
\unitlength 0.6mm
\begin{picture}(30,35)
\thinlines
\put(5,10){\dashbox{2}(20,10){$t$}}
\put(15,30){\vector(0,-1){10}}
\put(15,10){\vector(0,-1){10}}
\put(20,25){\makebox(0,0){$\rho$}}
\put(20,5){\makebox(0,0){$\sigma$}}
\end{picture}
\end{center}
\caption{An intertwiner $t\in\Hom(\rho,\sigma)$}
\label{trs}
\end{figure}
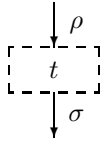

\begin{figure}[htb]
\begin{center}
\unitlength 0.6mm
\begin{picture}(42,35)
\thinlines
\put(5,10){\dashbox{2}(20,10){$t$}}
\put(15,30){\vector(0,-1){10}}
\put(15,10){\vector(0,-1){10}}
\put(35,30){\vector(0,-1){30}}
\put(20,25){\makebox(0,0){$\rho$}}
\put(20,5){\makebox(0,0){$\sigma$}}
\put(40,5){\makebox(0,0){$\tau$}}
\end{picture}
\end{center}
\caption{An intertwiner
$t\in\Hom(\rho,\sigma)\subset\Hom(\rho\tau,\sigma\tau)$}
\label{trst}
\end{figure}
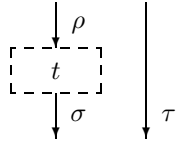

\begin{figure}[htb]
\begin{center}
\unitlength 0.6mm
\begin{picture}(40,35)
\thinlines
\put(15,10){\dashbox{2}(20,10){$t$}}
\put(25,30){\vector(0,-1){10}}
\put(25,10){\vector(0,-1){10}}
\put(5,30){\vector(0,-1){30}}
\put(30,25){\makebox(0,0){$\rho$}}
\put(30,5){\makebox(0,0){$\sigma$}}
\put(10,5){\makebox(0,0){$\mu$}}
\end{picture}
\end{center}
\caption{The intertwiner $\mu(t)\in\Hom(\mu\rho,\mu\sigma)$}
\label{mtrs}
\end{figure}
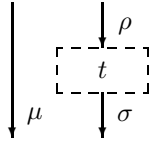

\begin{figure}[htb]
\begin{center}
\unitlength 0.6mm
\begin{picture}(170,40)
\thinlines
\put(8,5){\makebox(0,0){$\rho$}}
\put(15,20){\vector(0,-1){20}}
\put(25,20){\arc{20}{3.142}{0}}
\put(45,20){\arc{20}{0}{3.142}}
\put(55,40){\line(0,-1){20}}
\put(70,20){\makebox(0,0){$=$}}
\put(85,40){\vector(0,-1){40}}
\put(92,5){\makebox(0,0){$\rho$}}
\put(100,20){\makebox(0,0){$=$}}
\put(115,40){\line(0,-1){20}}
\put(125,20){\arc{20}{0}{3.142}}
\put(145,20){\arc{20}{3.142}{0}}
\put(155,20){\vector(0,-1){20}}
\put(162,5){\makebox(0,0){$\rho$}}
\end{picture}
\end{center}
\caption{A topological invariance}
\label{isoinv1}
\end{figure}
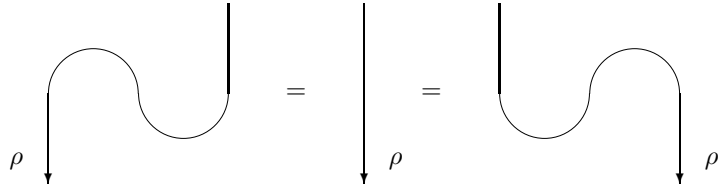

\begin{figure}[htb]
\begin{center}
\unitlength 0.5mm
\begin{picture}(205,30)
\thinlines
\put(10,0){\arc{20}{3.142}{0}}
\put(10,25){\arc{20}{0}{3.142}}
\put(20,25){\vector(0,1){0}}
\put(0,0){\vector(0,-1){0}}
\put(40,12.5){\makebox(0,0){,}}
\put(80,0){\arc{20}{3.142}{0}}
\put(80,25){\arc{20}{0}{3.142}}
\put(70,25){\vector(0,1){0}}
\put(90,0){\vector(0,-1){0}}
\put(60,25){\vector(0,-1){25}}
\put(110,12.5){\makebox(0,0){,}}
\put(160,0){\arc{20}{3.142}{0}}
\put(160,25){\arc{20}{0}{3.142}}
\put(170,25){\vector(0,1){0}}
\put(150,0){\vector(0,-1){0}}
\put(140,0){\vector(0,1){25}}
\put(130,25){\vector(0,-1){25}}
\put(190,12.5){\makebox(0,0){,}}
\put(200,12.5){\makebox(0,0){$\ldots$}}
\end{picture}
\end{center}
\caption{Jones projections in the tunnel}
\label{Jproj}
\end{figure}
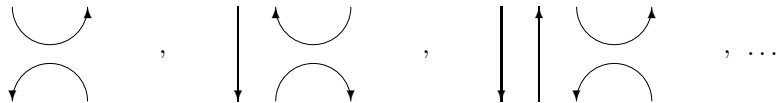

Sometimes it is useful to use graphical expressions for
formulae involving intertwiners. Roughly speaking,
for an intertwiner $t\in\Hom(\rho,\sigma)$ we draw
a picture as in \fig{trs}, i.e.\ we represent
morphisms by oriented ``wires'' and intertwiners by
boxes.
%
Reversing an arrow means replacing a label $\rho$ by
its conjugate $\co\rho$, and taking adjoints then
corresponds to vertical reflection
of the picture together with reversing all arrows.
As $\Hom(\rho,\sigma)\subset\Hom(\rho\tau,\sigma\tau)$ we
are allowed to add or remove straight wires on the right,
i.e.\ we are free to pass from \fig{trs}
to \fig{trst}.
%
On the other hand, the intertwiner $\mu(t)$ is in
$\Hom(\mu\rho,\mu\sigma)$ and is represented graphically
as in \fig{mtrs}.
%
With the convention that the identity morphism (of some factor)
is labelled by ``the invisible wire'', the isometries
$r_\rho$, ${\co r}_\rho$ and $r_\rho^*$, ${\co r}_\rho^*$
are represented as caps and cups, respectively, with different
orientations of the wire labelled by $\rho$.
Then, with certain normalization procedures taken care of
in \cite{BEK1} (where the graphical framework is worked
out in full detail -- but see also \cite{MS1,KR,Wi,FK,FG1,J2}),
the relations of \erf{riso} become
topological moves as in \fig{isoinv1}.
%
The minimal conditional expectation is obtained as
follows. First, the map
$\phi_\rho : B\rightarrow A$, $b\mapsto r_\rho^*\co\rho(b) r_\rho$,
is the unique standard left inverse for $\rho$
(as $\phi_\rho\circ\rho=\id_A$) and then
$E_\rho=\rho\circ\phi_\rho : B\rightarrow\rho(A)$
is the minimal conditional expectation for the
subfactor $\rho(A)\subset B$.
In the graphical framework, Jones projections in the
tunnel which were translates of $v_1v_1^*$ appear as
in \fig{Jproj}.
%
The Pimsner-Popa bound in the Kosaki-Jones index
is realized by such Jones projections so that the constant
$d_\rho$ in \erf{riso} is identified with
$[B:\rho(A)]^{1/2}$, the square root of the Jones index.

Returning to our original subfactor $N\subset M$
with inclusion homomorphism $\iota:N\hookrightarrow M$,
$\can=\iota\co\iota$, $\canr=\co\iota\iota$, where
$\co\iota$ is a conjugate for $\iota$, we have isometries
$w\equiv r_\iota\in\Hom(\id_N,\canr)$ and
$v\equiv {\co r}_\iota\in\Hom(\id_N,\can)$ satisfying
the consistency relations
$w^*v=w^*\can(v)=d_\iota^{-1}\bfe$ with
$d_\iota^2=[M:N]=d_\can=d_\canr$.
Note that we have pointwise equality
$M=Nv$ as $m=[M:N]^{1/2}w^*\can(m)v$ where
$[M:N]^{1/2}w^*\can(m)\in N$, $m\in M$, which means
that $v$ is a basis element for $M$ as an $N$-module.
The previous characterization of conjugates can be
used to characterize which endomorphisms arise as
a canonical endomorphism.

If $\can\in\End(M)$ where $M$ is an infinite factor,
then $\can$ is a canonical endomorphism of some
subfactor $N\subset M$ if and only if there exist
isometries $v\in\Hom(\id_M,\can)$ and
$w\in\Hom(\can,\can^2)$ such that
\begin{eqnarray}
w^*\can(w)=ww^* \,, \qquad \can(w)w=w^2 \,,\label{qsys1}\\
v^*w=w^*\can(v)=d^{-1}\bfe \,,\quad d>0 \,.\label{qsys2}
\end{eqnarray}
Note that if $v={\co r}_\iota$ and $w=r_\iota$ as before,
then $w\in\Hom(\id_N,\canr)\subset\Hom(\can,\can^2)$.
Conversely, if \erf{qsys1} and \erf{qsys2} hold then we
can define
$N=\{ x\in M : wx=\can(x)w \,,\,\, wx^*=\can(x^*)w\}$,
and then $E:M\rightarrow N$ defined by
$E(x)=w^*\can(x)w$, $x\in M$, is a conditional expectation.

\section{Subfactors arising from loop groups}

We now turn to the actual algebras which we will use to
describe our modular invariants, arising from loop groups.
The loop group $\LSUn$ consists of smooth maps
$f:S^1\rightarrow\SUn$, the product being pointwise multiplication.
The representations of interest will be projective representations
of $\LSUn$ which extend to positive energy representations of
$\LSUn\rtimes\mathrm{Rot}(S^1)$ where the rotation group acts
on the maps of $S^1$ in a natural way so that the ``Hamiltonian''
or infinitesimal generator $L_0$ is positive.
The ones of particular interest, the irreducible
unitary positive
energy representations are classified as follows.
First there is a level $k$, a positive integer describing
a cocycle because we are dealing with projective representations.
The projective representation restricts to a genuine
irreducible representation of the constant loops, identified
with $\SUn$ itself, the multiplier becomes irrelevant now
since we are dealing with simply connected groups.
In order to obtain positive energy, only finitely many
irreducible representations are admissible,
namely the vertices of (i.e.\ integrable weights in)
the Weyl alcove $\cA^{(n,k)}$.
The (adjacency matrices of the) graphs $N_\la$,
such as $N_\Box$ itself, describe the fusion of
positive energy representations.

Restricting to loops concentrated on an interval
$I\subset S^1$ (proper, i.e.\ $I\neq S^1$ and non-empty),
the corresponding subgroup denoted by
\[
\LISUn=\{f\in\LSUn : f(z)=1 \,,\,\, z\notin I \} \,,
\]
one finds that in each positive energy representation
$\pi_\la$ the sets of operators $\pi_\la(\LISUn)$
and $\pi_\la(\LIcSUn)$ commute where $I^\rmc$ is the
complementary interval of $I$, again using that
$\SUn$ is simply connected. In turn we obtain
a subfactor
\begin{equation}\label{laincl}
\pi_\la(\LISUn)''\subset\pi_\la(\LIcSUn)'\,,
\end{equation}
involving hyperfinite \typeiii$_1$ factors (see \cite{W}).
In the vacuum representation, labelled by $\la=0$,
we have Haag duality in that the inclusion collapses
to a single factor $N(I)=N(I)$, and in general we
obtain a subfactor. The level 1 representations
of $\LSUn$ are realized through the Fock state of the
Hardy space projection $P$ on $L^2(S^1;\bbC^n)$.
Since $[f,P]$ is Hilbert-Schmidt for $f\in\LSUn$
acting naturally on $L^2(S^1;\bbC^n)$, we have that
$\LSUn$ is implemented in the corresponding Fock space
giving a positive energy representation.

The vacuum representation $\pi_0$ gives a clear geometric
picture of the Tomita-Takesaki modular group action $\sigma$
and modular conjugation $J$ on the Fock vacuum vector,
cyclic and separating for say $\pi_0(\LISUn)''$ for $I$ being
the upper half circle. The Tomita-Takesaki modular group is
induced by the second quantization of the geometric action
of $\mathit{SU}(1,1)$ on $S^1$, which is seen to be ergodic.
Consequently the algebra must be \typeiii$_1$ as the
action $\sigma$ is never ergodic otherwise
(see e.g.\ \cite[Cor.\ 1.10.8]{Bg}).
Similarly the conjugation $J$ acts by flipping the circle,
taking $I$ into the complementary interval:
\[ \pi_0(\LISUn)''=J\pi_0(\LISUn)'J=\pi_0(\LIcSUn)' \,, \]
so that Haag duality holds in the vacuum representation
and is a consequence of Tomita-Takesaki theory.
More generally the inclusion \erf{laincl} can be
read as providing an endomorphism $\la$ (by abuse of
notation denoted by the same symbol as the label)
of the local algebra $N(I)$ such that \erf{laincl}
is isomorphic to $\la(N(I))\subset N(I)$.
By A.\ Wassermann's work \cite{W} we obtain this way
a system of endomorphisms $\sys=\{\la\}$, the morphisms
being labelled by the Weyl alcove $\cA^{(n,k)}$, which is
closed under sector fusion, and the fusion coefficients
$N_{\la,\mu}^\nu$ match exactly the loop group fusion.
Similar results have been obtained for minimal models \cite{Loke}
and (partially) ${\mathit{LSpin}}(2n)$ models \cite{Tol}.
(That the DHR morphisms of net of a conformal field theory
model obey exactly the Verlinde fusion rules from the
conformal character transformations was
conjectured in \cite{FG2}. Proofs for special cases
can be found in \cite{W,Loke,Tol,B1,B2}.
Antony Wassermann has informed us that he has computed
fusion for all simple, simply connected loop groups;
and with Toledano-Laredo all but E$_8$ using a variant of the
Dotsenko-Fateev differential equation considered in his thesis.)
If we take the relative commutants for the tunnel
\[
\cdots \subset \la\co\la\la(N) \subset \la\co\la(N)
\subset \la(N) \subset N \,,
\]
we are decomposing products $\la\co\la\la\cdots$ into
irreducible components obtaining the same relative commutants
as for the Jones-Wenzl \typeii$_1$ $\SUn_k$ subfactors.
More precisely, if $\la$ is the fundamental representation
$\Box$ and $A\subset B$ is the hyperfinite \typeii$_1$
subfactor
\begin{equation}\label{Heckesub}
\{g_i : i=1,2,3,\ldots\}''\subset\{g_i : i=0,1,2,\ldots\}''\,,
\end{equation}
where the $g_i$'s are the Hecke algebra generators
obtained as explained below,
then (using Popa \cite{P}) the loop group subfactor
$\la(N)\subset N$ is isomorphic to
$A\otimes N \subset B\otimes N$.

The statistical mechanical models of \cite{DJMO} are
generalizations of the Ising model. The configuration
space of the Ising model, distributions of
symbols ``$+$'' and ``$-$'' on the vertices of the
square lattice $\bbZ^2$, can be thought of edges on the
Dynkin diagram A$_3$ on the edges of a square lattice,
where the end vertices are labelled by ``$+$'' and ``$-$''.
This model can be generalized by replacing A$_3$ by
other graphs $\Gamma$ such as other Dynkin diagrams or
indeed the Weyl alcove $\cA^{(n,k)}$.
A configuration is then a distribution of the edges
of $\Gamma$ over $\bbZ^2$, and associated to each local
configuration is a Boltzmann weight satisfying the
Yang-Baxter equation. The justification of $\SUn$ models
is as follows. By Weyl duality, the representation of
the permutation group on $\bigotimes\Mat(n)$ is the
fixed point algebra of the product action of $\SUn$.
Deforming this, there is a representation of the Hecke
algebra in $\bigotimes\Mat(n)$, whose commutant is
a representation of a deformation of $\SUn$,
the quantum group $\SUn_q$ \cite{Ji}.
The Boltzmann weights lie in this algebra
representation, and at critically reduce to the
natural braid generators $g_i$, so that the
Yang-Baxter equation satisfied by the Boltzmann weights
reduces to the braid relation
$g_i g_{i+1} g_i = g_{i+1} g_i g_{i+1}$.
When $q=\E^{2\pi\I/(n+k)}$ is a root of unity,
the irreducible representations of the corresponding
Hecke algebra are labelled precisely by $\cA^{(n,k)}$.

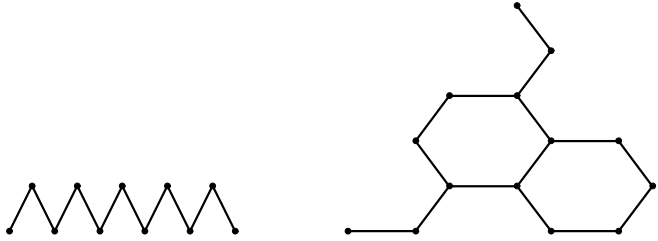
\begin{figure}[htb]
\begin{center}
\unitlength 0.15mm
\begin{picture}(600,230)
\thinlines
\multiput(0,0)(40,0){6}{\circle*{5}}
\multiput(20,40)(40,0){5}{\circle*{5}}
\thicklines
\multiput(0,0)(40,0){5}{\line(1,2){20}}
\multiput(20,40)(40,0){5}{\line(1,-2){20}}
\thinlines
\multiput(300,0)(60,80){2}{\circle*{5}}
\multiput(390,120)(60,80){2}{\circle*{5}}
\multiput(360,0)(30,40){2}{\circle*{5}}
\multiput(450,120)(30,40){2}{\circle*{5}}
\multiput(450,40)(30,40){2}{\circle*{5}}
\multiput(480,0)(60,80){2}{\circle*{5}}
\multiput(540,0)(30,40){2}{\circle*{5}}
\thicklines
\path(300,0)(360,0)(390,40)(360,80)(390,120)
(450,120)(480,160)(450,200)
\path(390,40)(450,40)(480,80)(450,120)
\path(450,40)(480,0)(540,0)(570,40)(540,80)(480,80)
\end{picture}
\end{center}
\caption{Colour zero-one part of the graphs in {\protect\fig{agraphs}}}
\label{agraphs01}
\end{figure}

The graph $\Gamma$ generates a von Neumann algebra by
considering larger and larger matrices generated by
larger and larger partition functions.
A subfactor can be obtained by the adjoint action,
placing the initial Boltzmann weights on the boundary.
For $\SUn_q$ subfactors this amounts to \erf{Heckesub}
because of the braid relations
$\Ad(g_1 g_2 \cdots)(g_i)=g_{i+1}$.
The principal graph of these inclusions is not the
entire graph $N_\Box$ (corresponding to $\cA^{(n,k)}$
as in \fig{agraphs})
but merely its zero-one part (the first two colours)
as in \fig{agraphs01}.
Nevertheless the entire graphs do have a meaning in subfactor
theory simply as graphs encoding the fusion rules of associated
systems of bimodules or sectors. Moreover, the center of $\SUn$,
namely $\bbZ_n$, induces an action on these subfactors and one may
construct crossed product or orbifold subfactors
$A^{\bbZ_n}\subset B^{\bbZ_n}$ \cite{EK1,X1} which will in turn
produce ``orbifold'' sector systems and graphs. As such graphs
have been noticed to label certain modular invariants, this can be
seen as a first indication that there is a relation between
modular invariant partition functions and subfactors. Another
strong indication is the special role of the Dynkin diagrams
$\rmD_{\mathrm{odd}}$ and E$_7$: In the classification of $\SUz_k$
modular invariants \cite{CIZ1,CIZ2,Kt}, the Dynkin diagrams A,
$\rmD_{\mathrm{even}}$, E$_6$ and E$_8$ label the \typei\
invariants whereas the invariants labelled by
$\rmD_{\mathrm{odd}}$ and E$_7$ are \typeii, i.e.\ involve a
non-trivial ``twist''. In subfactor theory it turned out that it
is precisely the diagrams A, $\rmD_{\mathrm{even}}$, E$_6$ and
E$_8$ which appear as principal graphs whereas
$\rmD_{\mathrm{odd}}$ and E$_7$ are not allowed (see \cite{Kaw}
and references therein).

\section{Braiding, $\a$-induction, and all that}

The geometry on the circle together with Haag duality
in the vacuum induces a braiding on the endomorphisms.
The endomorphisms $\la$ appearing above can be thought
of as being defined on a global algebra $\cN$ generated
by the $N(J)$'s where $J$ varies in the proper intervals
on $S^1$, neither touching nor containing a fixed
distinguished ``point at infinity'' $\zeta\in S^1$.
Then $\la$ will be localized on $I$ in the sense that
$\la(a)=a$ whenever $a\in N(J)$ with $J\cap I=\emptyset$,
and transportable in the sense that for each interval
$J$ there is a unitary $u\in\cN$ such that
$\Ad(u)\circ\la$ is localized in $J$.
Then if $\la$ and $\mu$ are localized on disjoint
intervals then they commute: $\la\mu=\mu\la$.
If however $\la$ and $\mu$ are localized in the same interval
$I$, then we may choose a relatively disjoint interval
$J$ (whose closure does not contain $\zeta$ as well)
and a unitary $u$ such that $\Ad(u)\circ\mu$ is
localized in $J$. Then $\la$ and $\Ad(u)\circ\mu$
commute and in turn $\eps_u(\la,\mu)=u^*\la(u)$ is a unitary
intertwining $\la\mu$ and $\mu\la$.
It turns out that this unitary is entirely independent
on the choice of $J$ and $u$, except that it may depend
on the choice of $J$ lying in the left or right
connected complement of $I$ with respect to the
point at infinity $\zeta$.
(See e.g.\ \cite{H,FRS1,FRS2,BE1} for more detailed
discussions of such matters.)
Therefore we have in fact only two ``statistics''
or braiding operators $\eps^+(\la,\mu)$ and
$\eps^-(\la,\mu)$, according to this choice.
Indeed we have $\eps^-(\la,\mu)=\eps^+(\mu,\la)^*$,
but $\eps^+(\la,\mu)$ and $\eps^-(\la,\mu)$ can be different.
The statistics operators obey a couple of consistency
equations which are called braiding fusion relations:
Whenever $t\in\Hom(\la,\mu\nu)$ one has
\[
\begin{array}{rl}
\rho(t) \eps^\pm (\lambda,\rho)
&= \eps^\pm(\mu,\rho) \mu(\eps^\pm(\nu,\rho))t \,, \\[.4em]
t \eps^\pm(\rho,\lambda)
&= \mu (\eps^\pm(\rho,\nu)) \eps^\pm(\rho,\mu) \rho(t) \,.
\end{array}
\]
This in turn implies the braid relation (or ``Yang-Baxter equation'')
\[
\rho(\eps^\pm(\lambda,\mu)) \eps^\pm(\lambda,\rho)
\lambda(\eps^\pm(\mu,\rho)) = \eps^\pm(\mu,\rho)
\mu(\eps^\pm(\lambda,\rho)) \eps^\pm(\lambda,\mu) \,.
\]
These equations turn our system $\sys$ of endomorphisms
into a ``braided $C^*$-tensor category'' (cf.\ \cite{DR}).

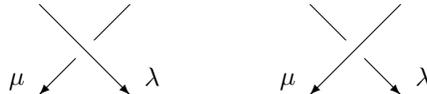
\begin{figure}[htb]
\begin{center}
\unitlength 0.6mm
\begin{picture}(100,25)
\thinlines
\put(10,20){\vector(1,-1){20}}
\put(30,20){\line(-1,-1){8}}
\put(18,8){\vector(-1,-1){8}}
\put(70,20){\line(1,-1){8}}
\put(82,8){\vector(1,-1){8}}
\put(90,20){\vector(-1,-1){20}}
\put(5,3){\makebox(0,0){$\mu$}}
\put(35,4){\makebox(0,0){$\la$}}
\put(65,3){\makebox(0,0){$\mu$}}
\put(95,4){\makebox(0,0){$\la$}}
\end{picture}
\end{center}
\caption{Braiding operators $\eps^+(\la,\mu)$ and
$\eps^-(\la,\mu)$ as over- and undercrossings}
\label{epslm}
\end{figure}

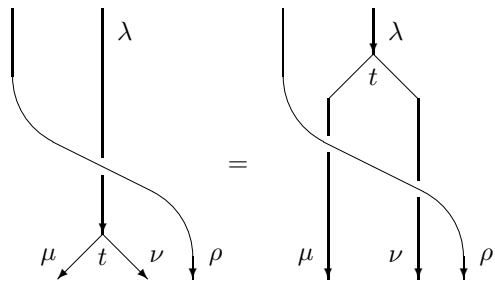
\begin{figure}[htb]
\begin{center}
\unitlength 0.6mm
\begin{picture}(120,65)
\thinlines
\put(26.180,45){\arc{32.361}{2.034}{3.142}}
\put(10,60){\line(0,-1){15}}
\put(30,25){\line(-2,1){11.1}}   
\put(33.820,5){\arc{32.361}{5.176}{6.283}}
\put(50,5){\vector(0,-1){5}}
\put(30,25){\line(2,-1){11.1}}
\put(30,60){\line(0,-1){33}}
\put(30,23){\vector(0,-1){13}}
\put(30,10){\vector(1,-1){10}}
\put(30,10){\vector(-1,-1){10}}
\put(30,5){\makebox(0,0){$t$}}
\put(35,55){\makebox(0,0){$\la$}}
\put(18,5){\makebox(0,0){$\mu$}}
\put(42,5){\makebox(0,0){$\nu$}}
\put(55,5){\makebox(0,0){$\rho$}}
\put(60,25){\makebox(0,0){$=$}}
\put(86.180,45){\arc{32.361}{2.034}{3.142}}
\put(70,60){\line(0,-1){15}}
\put(90,25){\line(-2,1){11.1}}
\put(93.820,5){\arc{32.361}{5.176}{6.283}}
\put(110,5){\vector(0,-1){5}}
\put(90,25){\line(2,-1){11.1}}
\put(90,60){\vector(0,-1){10}}
\put(90,50){\line(1,-1){10}}
\put(90,50){\line(-1,-1){10}}
\put(90,45){\makebox(0,0){$t$}}
\put(80,40){\line(0,-1){8}}
\put(100,40){\line(0,-1){18}}
\put(80,10){\line(0,1){18}}
\put(100,10){\line(0,1){8}}
\put(80,10){\vector(0,-1){10}}
\put(100,10){\vector(0,-1){10}}
\put(95,55){\makebox(0,0){$\la$}}
\put(75,5){\makebox(0,0){$\mu$}}
\put(95,5){\makebox(0,0){$\nu$}}
\put(115,5){\makebox(0,0){$\rho$}}
\end{picture}
\end{center}
\caption{The second braiding fusion equation for over-crossings}
\label{wireBFE2}
\end{figure}

\begin{figure}[htb]
\begin{center}
\unitlength 0.6mm
\begin{picture}(140,65)
\thinlines
\put(15,60){\line(0,-1){15.858}}
\put(25,30){\vector(1,-1){30}}
\put(25,30){\line(-1,1){7.071}}
\put(35,60){\line(1,-1){17.071}}
\put(55,30){\line(0,1){5.858}}
\put(55,30){\line(0,-1){5.858}}
\put(47,12){\line(1,1){5.071}}
\put(43,8){\vector(-1,-1){8}}
\put(55,60){\line(-1,-1){8}}
\put(43,48){\line(-1,-1){16}}
\put(23,28){\line(-1,-1){5.071}}
\put(15,15.858){\vector(0,-1){15.858}}
\put(25,44.142){\arc{20}{2.356}{3.142}}
\put(45,35.858){\arc{20}{5.498}{0}}
\put(45,24.142){\arc{20}{0}{0.785}}
\put(25,15.858){\arc{20}{3.142}{3.927}}
\put(8,5){\makebox(0,0){$\rho$}}
\put(30,5){\makebox(0,0){$\mu$}}
\put(60,5){\makebox(0,0){$\la$}}
\put(70,30){\makebox(0,0){$=$}}
\put(85,60){\line(1,-1){37.071}}
\put(125,15.858){\vector(0,-1){15.858}}
\put(105,60){\line(-1,-1){8}}
\put(93,48){\line(-1,-1){5.071}}
\put(85,30){\line(0,1){5.858}}
\put(85,30){\line(0,-1){5.858}}
\put(95,10){\vector(1,-1){10}}
\put(95,10){\line(-1,1){7.071}}
\put(125,60){\line(0,-1){15.858}}
\put(117,32){\line(1,1){5.071}}
\put(113,28){\line(-1,-1){16}}
\put(93,8){\vector(-1,-1){8}}
\put(95,24.142){\arc{20}{2.356}{3.142}}
\put(115,15.858){\arc{20}{5.498}{0}}
\put(115,44.142){\arc{20}{0}{0.785}}
\put(95,35.858){\arc{20}{3.142}{3.927}}
\put(80,5){\makebox(0,0){$\rho$}}
\put(110,5){\makebox(0,0){$\mu$}}
\put(132,5){\makebox(0,0){$\la$}}
\end{picture}
\end{center}
\caption{The braid relation as a vertical Reidemeister move of type III}
\label{Reid3}
\end{figure}
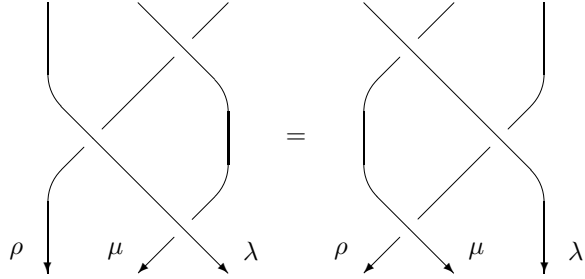

The braiding operators can be nicely incorporated in our
graphical intertwiner calculus. Namely, for $\eps^+(\la,\mu)$ and
$\eps^-(\la,\mu)$ we draw over- and undercrossings, respectively,
of wires $\la$ and $\mu$ as in \fig{epslm}.
%
Then the consistency relations are translated into
some kind of topological moves for the pictures,
as e.g.\ the second braiding fusion relation for
overcrossings is drawn graphically as in \fig{wireBFE2}
%
whereas the braid relation becomes a vertical Reidemeister
move of type III, presented in \fig{Reid3}.
%
We would like to obtain generators of the
modular group $\SLZ$ (up to normalization) from
the Hopf link and the twist, which is in fact
possible if and only if the braiding is
subject to a certain maximality condition,
called ``non-degeneracy'', basically stating that
$+$ and $-$ braiding operators are as different as
possible \cite{R1}.
For $\rho$ irreducible we find
$\eps^+(\rho,\co\rho)^*{\co r}_\rho=\om_\rho r_\rho$
for some scalar $\om_\rho\in\bbT$, thanks to
the uniqueness of isometries in the one-dimensional
$\Hom(\id,\co\rho\rho)$.

We will need net versions of canonical and dual canonical
endomorphisms to handle inclusions $N(I)\subset M(I)$,
where $N(I)$ are local, and which are standard in the
sense that there is a single vector $\Omega\in H$ being
cyclic and separating for all $M(I)$ on $H$ and all
$N(I)$ on a subspace $H_0\subset H$,
and such that there is a consistent family of conditional
expectations $E_I:M(I)\rightarrow N(I)$ preserving $\Omega$.
In this case for each interval $I_\circ$, there is an endomorphism
$\can$ of the global algebra $\cM$ associated to the net
$\{M(I)\}$ such that $\can|_{M(I)}$ is a canonical endomorphism
for $N(I)\subset M(I)$ whenever $I\supset I_\circ$.
Moreover, the restricted $\canr=\can|_\cN$ is localized
and transportable and we have
\begin{equation}\label{vacrestr}
\pi^0\simeq \pi_0\circ\can \,, \qquad
\pi^0|_\cN \simeq \pi_0\circ\canr \,,
\end{equation}
for $\pi^0$ denoting the defining representation of $\cM$ on $H$
and $\pi_0$ the representation of $\cN$ on $\overline{\cN\Omega}$.
There is an isometry $v$ intertwining the identity and $\can$,
and then we have $\cM=\cN v$, indeed $M(I)=N(I)v$ whenever
$I\supset I_\circ$. It is crucial to note that, though the net
$\{N(I)\}$ satisfies locality by assumption, the net $\{M(I)\}$
is not local in general.
In fact the latter is local if and only if the
chiral locality condition holds
\[ \eps^+(\canr,\canr)v^2=v^2 \,, \]
(see the original work \cite{LR} as an excellent mathematical
reference for these matters and \cite{RST} for a more physical
discussion of local extensions)
and locality of the extended net $\{M(I)\}$ is extremely
constraining, e.g.\ this automatically implies that all the
inclusions $N(I)\subset M(I)$ are irreducible, as shown
in \cite{BE1}.

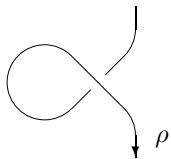
\begin{figure}[htb]
\begin{center}
\unitlength 0.5mm
\begin{picture}(43,45)
\thinlines
\put(10,20){\arc{20}{0.785}{5.498}}
\put(24.142,20){\line(-1,1){7.071}}
\put(22.142,18){\line(-1,-1){5.071}}
\put(24.142,20){\line(1,-1){7.071}}
\put(26.142,22){\line(1,1){5.071}}
\put(24.142,34.142){\arc{20}{0}{0.785}}
\put(24.142,5.858){\arc{20}{5.498}{6.283}}
\put(34.142,40){\line(0,-1){5.858}}
\put(34.142,5.858){\vector(0,-1){5.858}}
\put(41.142,5){\makebox(0,0){$\rho$}}
\end{picture}
\end{center}
\caption{Statistics phase $\om_\rho$ as a ``twist''}
\label{statph}
\end{figure}

\begin{figure}[htb]
\begin{center}
\unitlength 0.5mm
\begin{picture}(70,35)
\thinlines
\put(25,15){\arc{30}{5.742}{5.142}}
\put(45,15){\arc{30}{2.601}{2.001}}
\put(27,29.8){\vector(1,0){0}}
\put(43,29.8){\vector(-1,0){0}}
\put(5,25){\makebox(0,0){$\mu$}}
\put(65,25){\makebox(0,0){$\nu$}}
\end{picture}
\end{center}
\caption{Matrix element $Y_{\la,\mu}$ of Rehren's Y-matrix
         as a ``Hopf link''}
\label{Ymatrix}
\end{figure}
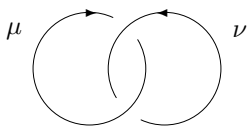

The statistics phase $\om_\rho$ ($\rho$ again irreducible)
can also be obtained using the left inverse
$\phi_\rho(\eps^+(\rho,\rho))=\om_\rho/d_\rho$.
Such formulae in algebraic quantum field theory
(see \cite{H} and references therein) predate subfactor theory.
Graphically $\om_\rho$ can be displayed as in \fig{statph}.
%
By a conformal spin and statistics theorem \cite{FG1,FRS2,GL2}
one can identify
\[ \om_\rho = \E^{2\pi\I h_\rho} \]
where $h_\rho$ is the lowest eigenvalue of the Hamiltonian
$L_0$ in the superselection sector $[\rho]$.
This will ensure that the statistics phase (and the
modular T-matrix) in our subfactor context coincide
with that in conformal field theory.
Now note that for $\mu,\nu$ irreducible the expression
$d_\mu d_\nu \phi_\mu (\eps(\nu,\mu)\eps(\mu,\nu))^*$
must be a scalar (as it is in $\Hom(\nu,\nu)$)
which we will denote by $Y_{\mu,\nu}$
and which is given graphically as in \fig{Ymatrix}.
%
In case we are dealing with a closed system $\sys$ of
braided endomorphisms it turns out \cite{R1,FG1,FRS2} that
\begin{equation}\label{Ymat}
Y_{\mu,\nu} =\sum_{\la\in\sys} \frac{\om_\mu\om_\nu}{\om_\la}
N_{\mu,\nu}^\la d_\la \,, \qquad \mu,\nu\in\sys\,.
\end{equation}
Normalizing the matrix $Y$ will yield the (modular) S-matrix.
Then from \erf{Ymat} it follows that if the $\om$'s and $N$'s
coincide then so does the modular matrix $S$ in the subfactor
context and that in conformal field theory.
Next define $z=\sum_{\la\in\sys} d_\la^2 \om_\la$.
If $z\neq 0$ put $c=4\arg(z)/\pi$, the central charge
which is defined modulo 8, and set
\[
S_{\la,\mu}= |z|^{-1} Y_{\la,\mu} \,, \qquad
T_{\la,\mu} = \E^{-\pi\I c/12} \om_\la \delta_{\la,\mu} \,.
\]
Then the matrices $S$ and $T$ obey the partial Verlinde
modular algebra
\[
TSTST=S \,, \quad CTC=T \,, \quad CSC=S \,,\quad
T^*T=\bfe \,,
\]
where $C$ is the conjugation matrix, i.e.\
$C_{\la,\mu}=\delta_{\la,\co\mu}$.
Moreover, the following conditions are equivalent
\cite{R1}:
\begin{itemize}
\item The braiding is non-degenerate, i.e.\
$\eps^+(\la,\mu)=\eps^-(\la,\mu)$ for all $\mu\in\sys$
only if $\la=\id$.
\item We have $|z|^2=w$ (recall that
$w=\sum_{\la\in\sys} d_\la^2$ is the global index
of the system $\sys$) and $S$ is invertible so that
$S$ and $T$ obey the full Verlinde modular algebra,
in particular $(ST)^3=S^2=C$, and $S$ diagonalizes the
fusion rules, i.e.\ the Verlinde formula of \erf{verlinde}
holds.
\end{itemize}

In our setting of a subfactor $N\subset M$ with a system
$\NXN$ of braided endomorphisms of $N$ we will show
how to induce endomorphisms of $M$. This method
corresponds to Mackey-induction in the group-subgroup
subfactor. The standard subfactor induction
$\la\mapsto\iota\la\co\iota$ will not be multiplicative
on sectors as e.g.\ the statistical dimension is
multiplied by $d_\canr$ --- so that in some sense we need
to divide out by $\canr$. This is achieved by the notion
of $\a$-induction which goes back to Longo and Rehren \cite{LR}
in the (nets of) subfactor setting, and it was studied
in \cite{BE1,BE2,BE3,BEK1,BEK2,BE4} and in
a similar framework (the relation is explained in \cite{X3})
also in \cite{X2}.

Note that $\can(v)\in\Hom(\canr,\canr^2)$ as
$v\in\Hom(\id,\can)$. Therefore the braiding fusion relations
can be applied to obtain
\[
\eps^\pm(\la,\canr)\la\can(v)\eps^\pm(\la,\canr)^*
=\canr(\eps^\pm(\la,\canr)^*)\can(v) \,,
\]
and as
\[
\eps^\pm(\la,\canr)\la\can(n)\eps^\pm(\la,\canr)^*
=\canr\la(n) \,,\qquad n\in N \,,
\]
we find by $M=Nv$ that
$\Ad(\eps^\pm(\la,\canr))\circ\la\can$ maps $M$ into
$\can(M)$, and so, in a more stream-lined notation,
\[
\a_\la^\pm = \co\iota^{\,-1} \circ \Ad
(\eps^\pm(\lambda,\canr)) \circ \lambda \circ \co\iota
\]
is a well-defined endomorphism of $M$ such that
$\a^\pm_\la(v)=\eps^\pm(\la,\canr)^* v$
and $\a^\pm_\la|_N=\la$.
The maps $\a^+$ and $\a^-$ are well-defined on sectors
and are multiplicative, additive and preserve conjugates:
\[
\a^\pm_{\la\mu}=\a^\pm_\la\a^\pm_\mu \,, \qquad
\overline{\a^\pm_\la}=\a^\pm_{\co\la} \,, \qquad
[\a^\pm_\nu] = [\a^\pm_{\nu_1}]\oplus[\a^\pm_{\nu_2}]
\]
for $[\nu]=[\nu_1]\oplus[\nu_2]$. In particular
the sectors $[\a^\pm_\la]$ commute; indeed
\[
\a^\pm_\mu \a^\pm_\la = \Ad(\eps^\pm(\la,\mu))
\circ \a^\pm_\la \a^\pm_\mu \,.
\]
In the restriction direction we write
\[
\sigma_\beta = \co\iota\circ\beta\circ\iota \equiv
\can\circ\beta|_N
\]
for $\beta$ an endomorphism of $M$.
Now $\sigma$ is additive on sectors and preserves conjugates
but it is not multiplicative (as e.g.\ $\sigma_\id=\canr$).
In general we have
\[
\langle \a^\pm_\la,\beta \rangle \le \langle
\la, \sigma_\beta \rangle \,,
\]
with equality in the case of chiral locality.
One has to be careful though for which endomorphisms $\beta$ of $M$
one is considering in this formula.
The inequality is true for any subsector of
$[\a^\pm_\la]$, $\la\in\NXN$.

To help compute such subsectors and their fusion rules,
one has the relation
\[
\langle \a^\pm_\la,\a^\pm_\mu \rangle \le \langle
\canr \la, \mu \rangle \,,
\]
again with equality in the case of chiral locality.
Note that we really have divided out by $\canr$, as
in the case of standard sector induction
$\la\mapsto\iota\la\co\iota$ we would have
$\langle\iota\la\co\iota,\iota\mu\co\iota\rangle
=\langle\canr^2\la,\mu\rangle$
by Frobenius reciprocity \cite{I2}.

We may also compare the two different ``chiral'' inductions
$\a^+$ and $\a^-$. Then $\a^+_\la=\a^-_\la$ is equivalent to
the monodromy being trivial, i.e.\
$\eps^+(\la,\canr)\eps^+(\canr,\la)=\bfe$.
Moreover, whenever chiral locality holds then
we even have that the chiral induced sectors
coincide, $[\a^+_\la]=[\a^-_\la]$, if and only
if the monodromy is trivial \cite{BE1}.
Nevertheless one has quite generally that
\[
\a^-_\mu \a^+_\la = \Ad(\eps^+(\la,\mu))\circ\a^+_\la \a^-_\mu \,,
\]
so that the sectors $[\a^-_\mu]$ and $[\a^+_\la]$ clearly
commute. Indeed even their subsectors commute and this
gives rise to a relative braiding symmetry between the
chiral induced sectors \cite{BE3}.

\section{Modular invariants, graphs and $\a$-induction}

The A-D-E classification of \cite{CIZ1,CIZ2,Kt} associates
a Dynkin diagram to each $\SUz$ modular invariant in such a
way that the multiplicities of the eigenvalues
$S_{1,\la}/S_{0,\la}$ of the associated graphs match the
diagonal entries $Z_{\la,\la}$ of the modular invariant.
Here $S$ is the modular S-matrix for $\SUz$ at level $k$,
and $\la$ just takes the values in the $\SUz_k$ spins
$\la\in\{0,1,2,....,k\}$.
For $\SUd$, Di Francesco and Zuber \cite{DZ1,DZ2,DiF}
sought graphs to describe
the modular invariants in an analogous way, guided partly by
the principle that the affine A-D-E diagrams correspond to the
finite subgroups of $\SUz$, and so began with fusion or McKay
graphs \cite{McK} of finite subgroups of $\SUd$ and sought
truncations with the correct eigenvalues --- a science
essentially based on trial and error.
Nevertheless they found a lot of interesting and puzzling
relations between graphs, fusion rules and coupling matrices,
giving the impetus to further research.
We illustrate our subfactor approach through analyzing one of
the exceptional $\SUz$ modular invariants which occurs at
level $k=10$. The modular invariant is
\begin{equation}\label{ZE6}
Z_{\mathrm{E}_6}=|\chi_0 + \chi_6|^2
+ |\chi_4 +\chi_{10}|^2 + |\chi_3+\chi_7|^2 \,.
\end{equation}
This invariant was labelled by the Dynkin diagram E$_6$ by
\cite{CIZ1} since the diagonal part
$\{\la : Z_{\la,\la}\neq 0\}$ of the invariant
is $\{0,3,4,6,7,10\}$ in this case are the Coxeter exponents
of E$_6$, i.e.\ the eigenvalues of the incidence
(or adjacency) matrix of E$_6$ are precisely
$\{ S_{1,\la}/S_{0,\la}=2\cos((\la+1)\pi/12) : \la=0,3,4,6,7,10\}$.
The E$_6$ modular invariant can be obtained from the conformal
embedding $\SUz_{10}\subset\SOf_1$, i.e.\ an inclusion of
$\SUz$ in $\SOf$ such that the level 1 positive energy
representations of $\LSOf$ decompose
into the level 10 representations of $\LSUz$, with finite
multiplicity. The loop group $\LSOf$ has three level 1
representations, the basic (b), vector (v) and spinor (s)
representation, with characters $\chi_\rmb$, $\chi_\rmv$
and $\chi_\rms$, respectively, decomposing as
\[
\chi_\rmb=\chi_0+\chi_6\,,\quad\chi_\rmv=\chi_4+\chi_{10}\,,
\quad \chi_\rms=\chi_3+\chi_7\,,
\]
on $\LSUz$. The diagonal invariant
$|\chi_\rmb|^2+|\chi_\rmv|^2+|\chi_\rms|^2$
of $\SOf_1$ then immediately produces the
exceptional E$_6$ invariant of $\SUz_{10}$ of \erf{ZE6}.
The positive energy representations of b, v, s of $\SOf_1$
satisfy the Ising fusion rules with b being the identity and
in particular fusion by s corresponds to the Dynkin diagram
A$_3$ as in \fig{A3}. Analogous to what we discussed
for $\LSUn$, they give rise to three endomorphisms in the
loop group subfactor setting of $\LSOf$ with the same
Ising fusion rules \cite{B2}.
The conformal embedding $\SUz_{10}\subset\SOf_1$ then gives
in the vacuum representation a net of subfactors
$\pi^0(\LISUz)''\subset \pi^0(\LISOf)''$ or $N(I)\subset M(I)$.
Over the net $\{N(I)\}$ we have a system of braided endomorphisms
$\{\la_j\}$ labelled by vertices (enumerated by $j=0,1,...,10$,
$\la_0=\id$) of the Dynkin diagram A$_{11}$, and braided endomorphisms
$\{\tau_\rmb=\id,\tau_\rmv,\tau_\rms\}$ over $\{M(I)\}$
corresponding to the vertices of A$_3$, where the graphs
A$_{11}$ and A$_3$ represent fusion by $\la_1$ and $\tau_\rms$.

\begin{figure}[htb]
\unitlength 0.65mm
\begin{center}
\begin{picture}(80,35)
\thinlines
\multiput(10,10)(15,0){5}{\circle*{2}}
\put(40,25){\circle*{2}}
\thicklines
\put(10,10){\line(1,0){60}}
\put(40,10){\line(0,1){15}}
\put(10,3){\makebox(0,0){$[\a^+_0]$}}
\put(25,3){\makebox(0,0){$[\a^+_1]$}}
\put(40,3){\makebox(0,0){$[\a^+_2]$}}
\put(55,3){\makebox(0,0){$[\a^+_9]$}}
\put(70,3){\makebox(0,0){$[\a^+_{10}]$}}
\put(40,32){\makebox(0,0){$[\varsigma]$}}
\end{picture}
\caption{$\SUz_{10}\subset\SOf_1$: Fusion graph
of $[\a^+_1]$ on the chiral system E$_6^+$}
\label{E6}
\end{center}
\end{figure}
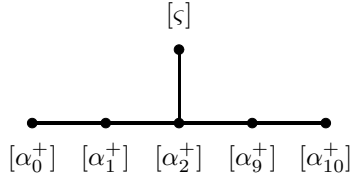

We can put the Ising A$_3$ system to one side for the time
being and focus on the A$_{11}$ system of $\{N(I)\}$.
Then (cf.\ \cite{LR,RST}) the dual canonical endomorphism
$\canr$ of $\cN$ is as a sector the sum $[\la_0]\oplus[\la_6]$
coming from the vacuum block --- this is basically \erf{vacrestr}.
We can thus perform first our $\a^+$-induction to obtain 11
endomorphisms $\{\a^+_j : j=0,1,...,10\}$
(we abbreviate $\a^+_j\equiv\a^+_{\la_j}$, but
after decomposition into irreducible sectors, we only find
six sectors $[\a^+_0]=[\id]$, $[\a^+_1]$, $[\a^+_2]$,
$[\a^+_9]$, $[\a^+_{10}]$ and $[\varsigma]$, the latter
appearing as a subsector of $[\a^+_3]$ which decomposes
as $[\a^+_3]=[\varsigma]\oplus[\a^+_9]$.
The graph E$_6$ appears as fusion graph of $[\a^+_1]$
as in \fig{E6}.
We now turn to our original braided system A$_3$ on $\{M(I)\}$.
Reading off from the blocks of the modular invariant
we find for $\sigma$-restriction:
\[
[\sigma_\rmb]=[\la_0]\oplus[\la_6]\equiv[\canr] \,,\quad
[\sigma_\rmv]=[\la_4]\oplus[\la_{10}] \,, \quad
[\sigma_\rms]=[\la_3]\oplus[\la_7]\,,
\]
because $\sigma$-restriction reflects the restriction
of (DHR) representations \cite{LR},
which is basically \erf{vacrestr} again.
Since the net of loop group factors $M(I)=\pi^0(\LISOf)''$
satisfies locality we have $\a\sigma$-reciprocity
\[
\langle \a^\pm_\la,\beta \rangle
= \langle \la, \sigma_\beta \rangle
\]
for $\la$ in the A$_{11}$ system and $\beta$ representing
any subsector of the induced system $\{[\a^+_j]\}$,
here E$_6$. Then as $\sigma$-restriction takes us
back into the A$_{11}$ system, the sectors
$[\tau_\rmb]$, $[\tau_\rmv]$, $[\tau_\rms]$ must lie
amongst the six E$_6$ sectors. They are identified
as $[\a^+_0]=[\tau_\rmb]$, $[\a^+_{10}]=[\tau_\rmv]$,
$[\varsigma]=[\tau_\rms]$, and indeed satisfy the
Ising fusion rules.

Other conformal inclusions and also simple current
extension invariants (often also called orbifold invariants)
can be handled similarly, the latter are realized by
so-called crossed product subfactors using the simple
current groups which represent the center $\bbZ_n$ of
$\SUn$ amongst the $\SUn_k$ fusion rules \cite{BE2,BE3}.
(See also Section \ref{realex}.)

So far we have only considered ``positive'' $\a^+$-induction,
arising from the braiding $\eps^+$. The same way we can use
the opposite braiding $\eps^-$, giving $\a^-$-induction
with ``negative'' chirality. In either case, say for the
conformal inclusion $\SUz_{10}\subset\SOf_1$, we have
two induced systems E$_6^+$ and E$_6^-$ of sectors on
$M(I)$, but at least they intersect on the Ising sectors
b, v, s of $\LSOf$ at level 1, symbolically:
$\rmE_6^+ \cap \rmE_6^- \supset \rmA_3$.
In fact they only coincide on these ``marked vertices'',
in the terminology of Di Francesco and Zuber, b, v, s
of E$_6^\pm$. Di Francesco and Zuber \cite{DZ1,DZ2,DiF}
had already empirically observed that the graphs which they
sought to describe the diagonal part of a given modular invariant
carried in the \typei\ case fusion rule algebras with certain
distinguished marked vertices forming fusion rule subalgebras
describing the extended fusion rules. This now finds
a clear explanation in terms of the game of induction and
restriction of sectors.

More generally, in the case of conformal embedding subfactors,
the following were shown to be equivalent \cite{BE3}:
\begin{itemize}
\item $\cV^+ \cap \cV^-=\cT$,
\item $Z_{\la,\mu}=\langle\a^+_\la,\a^-_\mu\rangle$,
\item The irreducible subsectors of $[\can]$ all lie
 in $\cV^\a=\cV^+\vee\cV^-$,
\item $\sum_{\beta\in\cV^\a} d_\beta^2=w$.
\end{itemize}
Here $\cV^\pm$ are the two chiral systems of induced irreducible
sectors, $\cT\subset\cV^\pm$ is the subsystem of neutral or
``ambichiral'' sectors, arising from either induction and
corresponding to the marked vertices, and finally $\cV^\a$
is the system of irreducible sectors generated by products
of the different chiral systems or equivalently obtained by
decomposing sectors $[\a^+_\la\a^-_\mu]$ into irreducibles.
The second condition gives a nice interpretation of the
modular invariant matrix $Z$ as counting the coupling of the
two chiral inductions. Note that it immediately produces
the upper bound of \erf{Zbound} because
the largest possible coupling occurs when $[\a^+_\la]$
and $[\a^-_\mu]$ both purely decompose into multiples
of one and the same irreducible sector and the multiplicities
are bounded by the statistical dimensions.
In fact the bound
$\langle\a^+_\la,\a^-_\mu\rangle\le d_\la d_\mu$
even holds for degenerate braidings
(i.e.\ with non-unitary S-matrices).
The third and hence all completeness properties could be
verified in a case by case analysis for e.g.\ all $\SUz$
and $\SUd$ conformal inclusion subfactors, and by a general
proof for all simple current extensions of $\SUn$ all levels
in \cite{BE3}. By adopting a graphical argument of
Ocneanu \cite{O} from the bimodule sectors to the sector
framework, the generating property was proven in \cite{BEK1}
to hold quite generally, provided the braiding is non-degenerate.
(And this is the case for $\SUn_k$ due to unitarity of the
S-matrix.)

A more careful analysis in \cite{BE4} using algebraic instead
of graphical techniques shows that the ``$\a$-global index''
$w_\a=\sum_{\beta\in\cV^\a} d_\beta^2$ is in fact given by
\[
w_\a = \frac w{\sum_{\la\,\,{\mathrm{deg}}} Z_{0,\la} d_\la}
\]
with summation over degenerate elements $\la$ for which
$\eps^+(\la,\mu)=\eps^-(\la,\mu)$ for all $\mu$.
Thus the generating property can hold even for some
degenerate systems. (An example is the conformal inclusion
subfactor $\SUz_{10}\subset\SOf_1$ if we start only with
the smaller system A$_{11}^{\mathrm{even}}$ of even spins.)
Moreover, the methods of \cite{BEK1,BEK2} allow us to handle
\typeii\ modular invariants as well as conformal embedding
and simple current invariants.

We now turn to the general framework of \cite{BEK1,BEK2,BE4}.
We take a subfactor $N\subset M$ and a system $\NXN\subset\End(N)$
of endomorphisms by which we mean a collection of irreducible
endomorphisms of finite statistical dimension, containing the
identity morphism and closed under conjugation and irreducible
decomposition of products. Then for $\iota:N\hookrightarrow M$
being the inclusion homomorphism and $\canr=\co\iota\iota$
and $\can=\iota\co\iota$ the dual canonical endomorphism
and canonical endomorphism, respectively, we assume that
$\canr$ lies in $\Sigma(\NXN)$, the set of morphisms
representing sector sums corresponding to the irreducibles
in $\NXN$ --- but make no assumption on $\can$. Moreover
we assume that the system $\NXN$ is braided.
We let $\MXM\subset\End(M)$ denote a system of endomorphisms
consisting of a choice of representative of each irreducible
subsector of sectors $[\iota\la\co\iota]$, $\la\in\NXN$.
We define $\MXMa\subset\MXM$ to be the subsystem of those
endomorphisms which are representatives of some
subsectors of $[\a^+_\la \a^-_\mu]$, $\la,\mu\in\NXN$.
(Note that by $\a^\pm_\la\iota=\iota\la$, any
subsector of $[\a^+_\la \a^-_\mu]$ will automatically
be a subsector of $[\iota\la\mu\co\iota]$ since
$[\can]$ contains the identity sector.)
Then we similarly define the chiral induced systems
as the subsystems $\MXMpm\subset\MXM$ of irreducible sectors
arising from positive/negative $\a^\pm$-induction,
and the neutral system $\MXMo=\MXMp\cap\MXMm$.
Their global indices, i.e.\ sums over squares of statistical
dimensions, are denoted by $w$, $w_\a$, $w_\pm$, and $w_0$
(it follows from the assumptions that $\NXN$ and $\MXM$
have the same global index $w$) and fulfill
$1 \le w_0 \le w_\pm \le w_\a \le w$.

{\em Defining} now a ``coupling matrix'' $Z$ by setting
\[
Z_{\la,\mu} = \langle \a^+_\la,\a^-_\mu \rangle \,,
\qquad \la,\mu\in\NXN \,,
\]
turns out to commute \cite{BEK1} with matrices $\Omega$ and $Y$,
where $\Omega_{\la,\mu}=\delta_{\la,\mu}\om_\la$ and
$Y$ is defined as in \erf{Ymat}.
(When the braiding is non-degenerate, we thus have a
physical modular invariant $Z$ which commutes with the
modular S- and T-matrices, being the normalized
matrices $Y$ and $\Omega$.)
Moreover, the relative
sizes of the various systems are encoded in $Z$,
namely we have \cite{BEK2}
\begin{equation}\label{chirglob}
w_+ = \frac w{\sum_{\la\in\NXN} d_\la Z_{\la,0}}
= \frac w{\sum_{\la\in\NXN} Z_{0,\la} d_\la} = w_-
\end{equation}
as well as \cite{BE4}
\[
w_\a = \frac w{\sum_{\la\in\NXNd} Z_{0,\la} d_\la}\,,
\qquad w_0=\frac{w_+^2}{w_\a} \,,
\]
where $\NXNd\subset\NXN$ denotes the subsystem of degenerate
morphisms (i.e.\ $\NXNd=\{\id\}$ in the non-degenerate case).
Here the equality
$\sum_\la d_\la Z_{\la,0} = \sum_\la Z_{0,\la} d_\la$
is due to the invariance $YZ=ZY$.

Although the original system $\NXN$ is braided, the
induced systems $\MXM$ or even $\MXMpm$ need not even be
commutative. Indeed if we complexify the fusion rules of
$\MXMpm$ to obtain finite-dimensional $C^*$-algebras
$\cZ^\pm$ we find \cite{BEK2}
(assuming non-degeneracy of the braiding)
\begin{equation}\label{cfuru}
\cZ^\pm \simeq \bigoplus_{\tau\in\MXMo}\bigoplus_{\la\in\NXN}
\Mat(b^\pm_{\tau,\la}) \,,
\end{equation}
where $b^\pm_{\tau,\la}=\langle\tau,\a^\pm_\la\rangle$
are the chiral branching coefficients for $\la\in\NXN$
and a neutral morphism $\tau\in\MXMo$ --- a marked vertex
of Di Francesco and Zuber.

In particular, the chiral systems are commutative only when
$b^\pm_{\tau,\la}\le 1$ for all $\tau,\la$. This explains
the non-commutativity discovered by Feng Xu \cite{X2}
with direct computations of some fusion rules for the
conformal embedding $\SUf_4\subset{\mathit{SU}}(15)_1$
(and which lead to conceptual problems in the partially
systematic approach to graphs from modular invariants of
\cite{PZ} based on certain assumptions)
and provides now a whole series of non-commutative chiral fusion
rules for $\SUn_n\subset{\mathit{SU}}(n^2-1)_1$, $n\ge 4$.
Moreover, by counting dimensions we find for the
cardinality $\#\MXMpm$ of $\MXMpm$ that
\[
\# {} \MXMpm = \sum_{\tau,\la} (b^\pm_{\tau,\la})^2 =
\tr (\tmat {b^\pm} b^\pm) \,.
\]
In the matrix algebra $\Mat(b^\pm_{\tau,\la})$, the
induced $[\a^\pm_\nu]$ is scalar, being $S_{\la,\nu}/S_{\la,0}$.
However, even in the degenerate case, the neutral elements
always possess a braiding (hence have commutative fusion)
arising as restriction of the relative braiding,
and this braiding is non-degenerate if the original
braiding on $\NXN$ is.
In that case we have ``extended '' S- and T-matrices
$S^\ext$ and $T^\ext$ from the neutral system,
and in $\Mat(b^\pm_{\tau,\la})$ a neutral sector $[\tau']$
($\tau'\in\MXMo$) acts as a central element since it
commutes with all subsectors of $[\a^\pm_\nu]$, as
$S^\ext_{\tau,\tau'}/S^\ext_{\tau,0}$.

Even if the chiral systems are commutative, the full system
$\MXMa=\MXMp\vee\MXMm$ may not be. Although $\MXMp$ and $\MXMm$
relatively commute (thanks to the relative braiding), it
may happen that ``mixed'' products of elements of $\MXMp$
and $\MXMm$ decompose into non-commuting irreducibles.
Indeed (cf.\ \erf{cfuru} for the chiral fusion rules),
if we complexify the fusion rules of $\MXMa=\MXM$ in the
non-degenerate case to obtain a finite-dimensional
$C^*$-algebra $\cZ$, then we find \cite{BEK1}
\[
\cZ \simeq \bigoplus_{\la,\mu\in\NXN} \Mat(Z_{\la,\mu}) \,.
\]
(This particular decomposition has also been claimed
by Ocneanu in his lectures \cite{O} in case of A-D-E
graphs and $\SUz$ modular invariants.)
Moreover, in the matrix algebra $\Mat(Z_{\la,\mu})$, the
induced $[\a^+_\nu]$ and $[\a^-_\nu]$ are scalars,
being $S_{\la,\nu}/S_{\la,0}$ and $S_{\mu,\nu}/S_{\mu,0}$,
respectively. We have seen in the case of chiral locality
(which holds e.g.\ for conformal embeddings) that we can
obtain graphs with spectrum corresponding to the diagonal
part of the modular invariant through the fusion graphs
of $[\a^\pm_\la]$ on $\MXMpm$.
In the general case where chiral locality may not necessarily
hold, we instead look at the action of $\MXM\supset\MXMpm$ on
the system $\MXN$ of $M$-$N$ sectors. Here the system
$\MXN$ is a choice of representatives of irreducible
subsectors of the sectors $[\iota\la]$, $\la\in\NXN$.
As $M$-$N$ sectors cannot be multiplied among themselves
there is no associated fusion rule algebra to decompose.
(Nevertheless, when chiral locality does holds, $\MXN$ can
be canonically identified with either $\MXMpm$ by
$\beta\mapsto\beta\circ\iota$, $\beta\in\MXMpm$.)
However, the left action of $\MXM$ on $\MXN$ defines a
representation $\varrho$ of the $M$-$M$ fusion rule algebra,
with matrix elements
$[\varrho([\beta])]_{\xi,\xi'}=\langle\xi,\beta\xi'\rangle$,
$\xi,\xi'\in\MXN$, and decomposes as
\cite{BEK1,BEK2}
\[
\varrho \simeq \bigoplus_{\la\in\NXN} \pi_{\la,\la} \,,
\]
where $\pi_{\la,\la}$ is the irreducible representation
corresponding to the matrix block $\Mat(Z_{\la,\la})$, so that
$\pi_{\la,\la}([\a^\pm_\nu])=S_{\la,\nu}/S_{\la,0}\bfe_{Z_{\la,\la}}$.
In particular the spectrum is determined by the diagonal part
of the modular invariant. Thus it is precisely this representation
$\varrho$ which provides an automatic connection between
the modular invariant and fusion graphs (e.g.\ the
representation matrix of some fundamental generator
$\Box$ corresponding to the left multiplication
of $[\a^\pm_{\Box}]$ on the $M$-$N$ sectors)
in such a way that (the multiplicities in) their spectra are
canonically given by the diagonal entries of the coupling matrix.
In fact, evaluation of $\varrho$ on the $[\a^\pm_\la]$'s
yields a ``nimrep'' of the original $N$-$N$ fusion rules,
i.e.\ a matrix representation where all the matrix entries
are non-negative integers.
Finally, by counting dimensions we see that $\#\MXN=\tr(Z)$.

We can illustrate this with the E$_7$ modular invariant
of $\SUz$:
\[
\begin{split}
Z_{{\mathrm{E}}_7} &= |\chi_0 + \chi_{16}|^2 +
|\chi_4 + \chi_{12}|^2 + |\chi_6 + \chi_{10}|^2 \\
& \qquad +|\chi_8|^2 + (\chi_2 + \chi_{14})\chi_8^* +
\chi_8(\chi_2 + \chi_{14})^* \,.
\end{split}
\]
Instead of simply extending to a diagonal invariant,
as in the E$_6$ case, we also insert a twist on the blocks.
This is an example of the setting of Moore and Seiberg \cite{MS2}
(see also Dijkgraaf and Verlinde \cite{DV}) that taking a
maximal extension of the ``chiral algebra'' $\cA\subset\cB$,
a modular invariant of $\cA$ is the restriction of some
permutation invariant
$Z^\ext_{\tau,\tau'}=\delta_{\tau,\omega(\tau')}$
where $\omega$ is a permutation of the sectors of the
extended theory $\cB$, defining an automorphism of
their fusion rules and preserving the extended vacuum
sector, $\omega(0)=0$.
The E$_7$ invariant is a twist of the D$_{10}$ invariant,
the latter we can realize from a subfactor with the dual canonical
endomorphism sector decomposing as $[\la_0]\oplus[\la_{16}]$,
being a simple current extension \cite{BE2}.
As shown in \cite{BEK2}, the E$_7$ modular invariant appears
for a subfactor with dual canonical endomorphism sector
$[\la_0]\oplus[\la_8]\oplus[\la_{16}]$.
For either invariant we find
$\tr(\tmat{b}^\pm b^\pm)=10$ ($=\tr(Z_{{\mathrm{D}}_{10}})$)
so that indeed in either case the fusion graph of the
generator $[\a^\pm_1]$ on $\MXMpm$ is D$_{10}$.
However, $\tr(Z_{{\mathrm{E}}_7})=7$ and the fusion graph
of $[\a^\pm_1]$ on $\MXN$ is E$_7$.

\section{Type II modular invariants,
extended fusion rule automorphisms, and all that}
\label{t2ext}

Now in our general setting we have
\[
Z_{\la,\mu} = \langle\a^+_\la,\a^-_\mu\rangle
= \sum_{\tau\in\MXMo} b^+_{\tau,\la} b^-_{\tau,\mu} \,,
\]
with chiral branching coefficients
$b^\pm_{\tau,\la}=\langle\tau,\a^\pm_\la\rangle$.
To write this in Moore-Seiberg form we would need
$b^-_{\tau,\la}=b^+_{\omega(\tau),\la}$ for
a permutation of the extended system, being identified
as the neutral system $\MXMo$,
so that
\[
Z_{\la,\mu} = \sum_{\tau\in\MXMo} b^+_{\tau,\la}
b^+_{\omega(\tau),\mu} \,.
\]
Note that by $\omega(0)=0$ and $b^\pm_{\tau,0}=\delta_{\tau,0}$
(do not worry that we denote both the original and
the extended ``vacuum'' i.e.\ identity morphism
by the same symbol ``$0$'')
we are automatically forced to have symmetric vacuum coupling
$Z_{\la,0}=Z_{0,\la}$.
To cover more general cases, which do occur as we shall see,
we should consider instead of one maximal extension
$\cA\subset\cB$ of the chiral algebra $\cA$, but two
different extensions $\cA\subset\cB_\pm$, yielding
different labelling sets of extended fusion rules so
that the extended modular invariant is
\[
Z^\ext_{\tau_+,\tau_-} = \delta_{\tau_+,\omega(\tau_-)} \,,
\]
where $\omega$ now is an isomorphism between the two sets
of extended fusion rules, still subject to $\omega(0)=0$.
Note that when we have two different labelling sets it
makes no sense to ask whether a coupling matrix is
symmetric or not.

When chiral locality does hold then
\[
b^\pm_{\beta,\la} = \langle\a^\pm_\la,\beta\rangle
=\langle\la,\sigma_\beta\rangle \,,
\]
whenever $\beta\in\MXMpm$. In particular, when
$\beta=\tau$ is neutral, i.e.\ lies in the intersection
$\MXMo=\MXMp\cap\MXMm$, then
\[
b^+_{\tau,\la}=b^-_{\tau,\la}\equiv b_{\tau,\la}\,,
\]
and we have a block decomposition or ``\typei''
modular invariant
\[
Z_{\la,\mu} = \sum_{\tau\in\MXMo} b_{\tau,\la}b_{\tau,\mu} \,.
\]
Permutation invariants can be classified as follows.
The following conditions are equivalent \cite{BEK2}:
\begin{itemize}
\item $Z_{\la,\mu}=\delta_{\la,\omega(\mu)}$ with
$\omega$ a permutation of $\NXN$ with $\omega(0)=0$
and defining a fusion rule automorphism,
\item $Z_{\la,0}=\delta_{\la,0}$,
\item $Z_{0,\la}=\delta_{\la,0}$,
\item $w_\pm=w$.
\end{itemize}
In this case the two inductions $\a^\pm$ are
isomorphisms (i.e.\ each $[\a^\pm_\la]$ is irreducible)
and $\omega=(\a^+)^{-1}\circ\a^-$. This result does
not rely on non-degeneracy of the braiding.

We would like to decompose a modular invariant into its
two parts, a \typei\ part together with a twist, and in
order to take care of heterotic vacuum coupling we will
need to implement such a twist by an isomorphism rather
than an automorphism. First we characterize chiral locality.
If chiral locality holds, i.e.\ $\eps^+(\canr,\canr)v^2=v^2$,
then
$Z_{\la,0}=\langle\a^+_\la,\a^-_0\rangle=\langle\a^+_\la,
\id\rangle=\langle\la,\canr\rangle$, and similarly
$Z_{0,\la}=\langle\la,\canr\rangle$. Indeed the following
conditions are equivalent \cite{BE4}:
\begin{itemize}
\item We have $Z_{\la,0}=\lan\canr,\la\ran$ for all $\la\in\NXN$.
\item We have $Z_{0,\la}=\lan\canr,\la\ran$ for all $\la\in\NXN$.
\item Chiral locality holds: $\eps^+(\canr,\canr)v^2=v^2$.
\end{itemize}
Thus chiral locality holds if and only if
\[
[\canr]=\bigoplus_{\la\in\NXN}\lan\canr,\la\ran[\la]
=\bigoplus_{\la\in\NXN}Z_{\la,0}[\la]
=\bigoplus_{\la\in\NXN}Z_{0,\la}[\la] \,.
\]
In general we define sectors
\[
[\canr_+]=\bigoplus_{\la\in\NXN}Z_{\la,0}[\la]\,,\qquad
[\canr_-]=\bigoplus_{\la\in\NXN}Z_{0,\la}[\la]\,.
\]
Note that
$d_{\canr_+}=\sum_\la d_\la Z_{\la,0}
=\sum_\la Z_{0,\la} d_\la=d_{\canr_-}$
(due to $(YZ)_{0,0}=(ZY)_{0,0}$) but in general
$[\canr_+]$ and $[\canr_-]$ may be different.
Using results on intermediate subfactors \cite{ILP}
it was shown in \cite{BE4} that, starting with
an arbitrary subfactor $N\subset M$ subject to
our assumptions, both $[\canr_+]$ and $[\canr_-]$
are dual canonical endomorphism sectors of $N$,
corresponding to intermediate subfactors
\[
N\subset M_\pm\subset M \,,
\]
and that $N\subset M_\pm$ satisfy chiral locality.
We then can form the $\tilde{\a}^\pm_\delta$-inductions
($\la\mapsto\tilde{\a}^\pm_{\delta;\la}$), on $N\subset M_\delta$,
$\delta=\pm$, and consider then the symmetric \typei\
modular invariants $Z^\pm$,
\[
Z^+_{\la,\mu}=\langle\tilde{\a}^+_{+;\la},
\tilde{\a}^-_{+;\mu}\rangle \,, \qquad
Z^-_{\la,\mu}=\langle\tilde{\a}^+_{-;\la},
\tilde{\a}^-_{-;\mu}\rangle \,.
\]
{}From the definition of $[\canr_\pm]$ we have
$Z^+_{\la,0}=\langle\canr_+,\la\rangle=Z^+_{0,\la}$
as $N\subset M_+$ satisfies chiral locality,
and so $Z^+_{\la,0}=Z^+_{0,\la}=Z_{\la,0}$
and similarly $Z^-_{\la,0}=Z^-_{0,\la}=Z_{0,\la}$.
For $Z=Z_{{\mathrm{E}}_7}$ the E$_7$ modular invariant
of $\SUz$, $Z^\pm$ will both be $Z_{{\mathrm{D}}_{10}}$
but it is possible as we shall see that $Z^+\neq Z^-$.

Next we argue that we can canonically identify
$\MXMp$ with $\MXMpp$ and $\MXMm$ with $\MXMmm$.
To do this it will be enough to find an injective
map $\MXMp\rightarrow\MXMpp$ (and $\MXMm\rightarrow\MXMmm$)
because the global indices $w^+$ are the same
thanks to \erf{chirglob} and $Z_{\la,0}=Z^+_{\la,0}$.
One can show that
$\Hom(\id,\a^\pm_\nu)=\Hom(\id,\tilde{\a}^\pm_{\pm;\nu})$
which in turn implies
$\Hom(\a^\pm_\la,\a^\pm_\mu)=
\Hom(\tilde{\a}^\pm_{\pm;\la},\tilde{\a}^\pm_{\pm;\mu})$.
In particular
$\Hom(\a^\pm_\la,\a^\pm_\mu)\subset M_\pm$.
We can then move from intertwiners to endomorphisms.
If $\beta\in\MXMp$ represents a subsector of $[\a^+_\la]$,
so that there is a $t\in M$ such that
$\a^+_\la(\cdot)=t\beta(\cdot)t^*+\ldots$, then
$tt^*\in\Hom(\a^\pm_\la,\a^\pm_\la)=
\Hom(\tilde{\a}^\pm_{\pm;\la},\tilde{\a}^\pm_{\pm;\la})$.
We can then construct an endomorphism
$\tilde{\beta}\in\End(M_+)$ representing a subsector of
$[\tilde{\a}^\pm_{\pm;\la}]$ and such that
$\beta|_{M_+}=\tilde{\beta}$.
In this way we construct bijections
$\vartheta_\pm : \MXMpm\rightarrow\MXMpmpm$
which \cite{BE4}:
\begin{itemize}
\item preserve chiral branching rules
 $\langle\beta,\a^\pm_\la\rangle=\langle\vartheta(\beta),
\tilde{\a}^\pm_{\pm;\la}\rangle$, $\beta\in\MXMpm$,
\item preserve chiral fusion rules,
\item and restrict to bijections of the neutral
 systems $\MXMo\rightarrow\MXMopm$.
\end{itemize}
This means that $\MXMo$ can be used (rather than $\MXMopm$)
to decompose the \typei\ coupling matrices
\[
Z^\pm_{\la,\mu} = \sum_{\tau\in\MXMo}
b^\pm_{\tau,\la} b^\pm_{\tau,\mu}
\]
with chiral branching coefficients
$b^\pm_\la=\langle\tau,\a^\pm_\la\rangle$,
$\tau\in\MXMo$, $\la\in\NXN$. If the two intermediate
subfactors happen to be identical, $M_+=M_-$
(so that the ``parent'' coupling matrices coincide,
$Z^+=Z^-$), then we can write
\[
Z_{\la,\mu} = \sum_{\tau\in\MXMo}
b^+_{\tau,\la} b^+_{\omega(\tau),\mu}
\]
for the (generically \typeii) coupling matrix $Z$.
Here the permutation
$\omega=\vartheta_+^{-1}\circ\vartheta_-$,
satisfying $\omega(0)=0$ clearly defines an automorphism
of the neutral fusion rules.

In general when $M_+\neq M_-$ we would write the
extended coupling matrix as
\[
Z^\ext_{\tau_+,\tau_-} = \delta_{\tau_-,\vartheta(\tau_+)} \,,
\]
where $\tau_\pm\in\MXMopm$ and
$\vartheta=\vartheta_-\circ\vartheta_+^{-1} :
\MXMop\rightarrow\MXMom$ is a bijection defining
an isomorphism of the chiral fusion rules.
We will illustrate that such heterotic situations do exist,
in fact examples are already provided by certain $\SOn_k$
current algebra models. We will deal with the simplest case
at level $k=1$ in Section \ref{hetero}.

What is the connection between the two chiral inductions
and the picture of left- and right-chiral algebras in
conformal field theory?
An appropriate notion of chiral algebras in the
setting of algebraic quantum field theory are
``chiral observables'' \cite{R5}, and one can show
that our coupling matrices describe in fact a Hilbert space
decomposition of the vacuum sector of a two-dimensional
quantum field theory upon restriction to the action of a
tensor product of left- and right-chiral observables \cite{R7}.
Suppose that our factor $N$ is obtained as a local factor
$N=N(I_\circ)$ of a quantum field theoretical net of factors
$\{N(I)\}$ indexed by proper intervals $I\subset \bbR$
on the real line, and that the system $\NXN$ is obtained
as restrictions of DHR-morphisms (cf.\ \cite{H}) to $N$.
This is in fact the case in our examples arising from
conformal field theory where the net is defined in terms
of local loop groups in the vacuum representation.
Taking two copies of such a net and placing the real axes
on the light cone, then this defines a local net $\{A(\cO)\}$,
indexed by double cones $\cO$ on two-dimensional
Minkowski space (cf.\ \cite{R5} for such constructions).
Given a subfactor $N\subset M$, determining in turn two
subfactors $N\subset M_\pm$ obeying chiral locality,
will provide two local nets of subfactors
$\{N(I)\subset M_\pm(I) \}$ as a local subfactor
basically encodes the entire information about the
net of subfactors \cite{LR}.
Arranging $M_+(I)$ and $M_-(J)$ on the two light
cone axes defines a local net of subfactors
$\{A(\cO)\subset A_\ext(\cO)\}$ in Minkowski space.
Rehren has recently proven \cite{R7} (see also \cite{BEK3}
for a different but less general derivation) that there is a
(\typeiii) factor $B$ such that we have an irreducible
inclusions $N\otimes N^\op\subset B$ such that the
dual canonical endomorphism $\Theta$ of the
inclusion $N\otimes N^\op \subset B$ decomposes as
\[
[\Theta] =\bigoplus_{\la,\mu\in\NXN}
Z_{\la,\mu} \, [\la \otimes \mu^\op] \,.
\]
(Here the superscript ``opp'' just denotes the opposite algebra,
i.e.\ $N^\op$ is $N$ as a linear space, with reversed
multiplication. There is a canonical way of identifying
$N(I)^\op$ with the CPT reflection of $N(I)$ \cite{GL1}
which is involved in the two-dimensional construction.)
Refining this result it has been shown \cite{BE4} that our
local extensions $M_\pm$ produce an intermediate subfactor
\[
N\otimes N^\op \subset M_+ \otimes M_-^\op \subset B
\]
such that moreover the dual canonical endomorphism
$\Theta_\ext$ of the inclusion $M_+\otimes M_-^\op \subset B$
decomposes as
\[
[\Theta_\ext] =\bigoplus_{\tau\in\MXMo}
[\vartheta_+(\tau) \otimes \vartheta_-(\tau)^\op] \,.
\]
The embedding $M_+\otimes M_-^\op \subset B$
gives rise to another net of subfactors
$\{A_\ext(\cO) \subset B(\cO)\}$, and a condition
which ensures that the net $\{B(\cO)\}$ obeys local
commutation relations can be established.
The existence of the local net was already proven
in \cite{R7}, and now the decomposition of $[\Theta_\ext]$
tells us that the chiral extensions $N(I)\subset M_+(I)$
and $N(I)\subset M_-(I)$ for left and right chiral nets
are indeed maximal (in the sense of \cite{R5}),
following from the fact that the coupling matrix for
$\{A_\ext(\cO) \subset B(\cO)\}$ is a bijection.
This shows that the inclusions $N\subset M_\pm$ should in
fact be regarded as the subfactor version of left- and
right maximal extensions of the chiral algebra.

\section{Heterotic examples}
\label{hetero}

Let us now consider the $\SOn$ loop group models at level 1,
where $n$ is a multiple of 16, $n=16\ell$, $\ell=1,2,3,...\,$.
These theories have four sectors, the basic ($0$), vector (v),
spinor (s) and conjugate spinor (c) module, corresponding to
highest weights $0$, $\Lambda_{(1)}$, $\Lambda_{(r-1)}$
and $\Lambda_{(r)}$, respectively; here
$r=n/2=8\ell$ is the rank of $\SOn$. The conformal
dimensions are given as $h_0=0$, $h_\rmv=1/2$,
$h_\rms=h_\rmc=\ell$, and the sectors obey
$\bbZ_2\times\bbZ_2$ fusion rules. The
Kac-Peterson matrices are given explicitly as
\begin{equation}\label{ST}
S = \frac 12 \left( \begin{array}{rrrr}
1 & 1 & 1 & 1 \\ 1 & 1 & -1 & -1 \\
1 & -1 & 1 & -1 \\ 1 & -1 & -1 & 1
\end{array} \right) , \qquad
T = \E^{-2\pi\I \ell/3} \left( \begin{array}{rrrr}
1 & 0 & 0 & 0 \\ 0 & -1 & 0 & 0 \\
0 & 0 & 1 & 0 \\ 0 & 0 & 0 & 1
\end{array}  \right) .
\end{equation}
It is easy to check that there are exactly six
modular invariants, $Z=\bfe$, $W$, $X_\rms$,
$X_\rmc$, $Q$, $\tmat Q$. Here
\[
W = \left( \begin{array}{rrrr}
1 & 0 & 0 & 0 \\ 0 & 1 & 0 & 0 \\
0 & 0 & 0 & 1 \\ 0 & 0 & 1 & 0
\end{array}  \right), \;\;
X_\rms = \left( \begin{array}{rrrr}
1 & 0 & 1 & 0 \\ 0 & 0 & 0 & 0 \\
1 & 0 & 1 & 0 \\ 0 & 0 & 0 & 0
\end{array}  \right), \;\;
Q = \left( \begin{array}{rrrr}
1 & 0 & 0 & 1 \\ 0 & 0 & 0 & 0 \\
1 & 0 & 0 & 1 \\ 0 & 0 & 0 & 0
\end{array}  \right),
\]
and $X_\rmc=W X_\rms W$.
(Note that $Q=X_\rms W$ and $\tmat Q=W X_\rms$.)
The matrix $Q$ and its transpose $\tmat Q$
are two examples for modular invariants with
non-symmetric vacuum coupling. Such ``heterotic''
invariants seem to be extremely rare and have not
enjoyed particular attention in the literature,
perhaps because they were erroneously dismissed
as being spurious in the sense that they would not
correspond to a physical partition function.
Examples for truly spurious modular invariants
were given in \cite{SY2,V,FSS} and found
to be ``coincidental'' linear combinations of
proper physical invariants. Note that although
there is a linear dependence here, namely
\[
\bfe - W - X_\rms - X_\rmc + Q + \tmat Q = 0 \,,
\]
we cannot express $Q$ (or $\tmat Q$) alone as a linear
combination of the four symmetric invariants.
This may serve as a first indication that $Q$ and
$\tmat Q$ are not spurious. We will now demonstrate
that they can be realized from subfactors.

The $\bbZ_2\times\bbZ_2$ fusion rules for these models
were proven in the DHR framework in \cite{B2}, and
together with the conformal spin and statistics theorem
\cite{FG1,FRS2,GL2} we conclude that there is a net of
\typeiii\ factors on $S^1$ with a system
$\{\id,\rho_\rmv,\rho_\rms,\rho_\rmc\}$ of localized
and transportable, hence braided endomorphisms, such that
the statistics S- and T-matrices are given by \erf{ST}.
Because the statistics phases are second roots of unity as
$\omega_\rmv=-1$ and $\omega_\rms=\omega_\rmc=1$,
we can by \cite{R3} choose the morphisms in the system
such that obey the $\bbZ_2\times\bbZ_2$ fusion rules even
by individual multiplication,
\[
\rho_\rmv^2=\rho_\rms^2=\rho_\rmc^2=\id\,, \qquad
\rho_\rmv \rho_\rms = \rho_\rms \rho_\rmv = \rho_\rmc \,.
\]
This is enough to proceed with the DHR construction of
the field net \cite{DHR2}, as already carried out
similarly for simple current extensions with cyclic
groups in \cite{BE2,BE3}. In fact, all we need to do here
is to pick a single local factor $N=N(I)$ such that the
interval $I\subset S^1$ contains the localization region
of the morphisms, and then we construct the cross product
subfactor $N\subset N\rtimes(\bbZ_2\times\bbZ_2)$.
Then the corresponding dual canonical endomorphism
$\canr$ decomposes as a sector as
\[
[\canr] = [\id] \oplus [\rho_\rmv] \oplus
[\rho_\rms] \oplus [\rho_\rmc] \,.
\]
Checking
$\lan\iota\la,\iota\mu\ran=\lan\canr\la,\mu\ran=1$
for $\la,\mu=\id,\rho_\rmv,\rho_\rms,\rho_\rmc$,
we find that there is only a single $M$-$N$ sector,
namely $[\iota]$. From $\tr Z=\#\MXN$ we conclude
that the modular invariant coupling matrix $Z$
arising from this subfactor must fulfill $\tr Z=1$.
This leaves only the possibility that $Z$ is $Q$
or $\tmat Q$. We may and do assume that $Z=Q$,
otherwise we exchange braiding and opposite braiding.
It is easy to determine the intermediate subfactors
$N\subset M_\pm \subset M$. Namely, we have
$M_+=N\rtimes_{\rho_\rms} \bbZ_2$ and
$M_-=N\rtimes_{\rho_\rmc} \bbZ_2$ with dual
canonical endomorphism sectors
$[\canr_+]=[\id]\oplus[\rho_\rms]$ and
$[\canr_-]=[\id]\oplus[\rho_\rmc]$, respectively.
That both extensions are local also follows from 
$\omega_\rms=\omega_\rmc=1$. We therefore
find $Z^+=X_\rms$ and $Z^-=X_\rmc$. Finally, the
permutation invariant $W$ is obtained from the
non-local extension $M_\rmv=N\rtimes_{\rho_\rmv}\bbZ_2$.

\section{Realization of modular invariants from subfactors}
\label{realex}

In our general setting, we have the following situation:
For a given \typeiii\ von Neumann factor $N$ equipped with
a braided system of endomorphism $\NXN$, any embedding
$N\subset M$ of $N$ in a larger factor $M$ which is
compatible with the system $\NXN$ (in the sense that
the dual canonical endomorphism decomposes in $\NXN$)
defines a coupling matrix $Z$ through $\a$-induction.
This matrix $Z$ commutes with the matrices $Y$ and $\Omega$
arising from the braiding and in turn is a
``modular invariant mass matrix'' whenever the braiding
is non-degenerate. Suppose we start with a system
corresponding to the RCFT data of $\SUn_k$.
Then the following question is natural, but difficult to answer:
\begin{quote}
\textbf{Can any physical modular invariant be
realized from some subfactor $N\subset M$?}
\end{quote}
\noindent The first problem with this question is that
one needs to specify what the term ``physical'' means.
Quite often in the literature, any modular invariant matrix
(i.e.\ $ZS=SZ$, $ZT=TZ$) subject to the constraint that
all entries are non-negative integers and with normalization
$Z_{0,0}=1$ is called a physical invariant.
Well, with this interpretation of ``physical'' the answer
to the question is clearly negative. Namely, our
general theory says that there is always some associate
extended theory carrying another representation of
the modular group $\SLZ$ which is compatible with
the chiral branching rules.
As mentioned above, it is however known \cite{SY2,V,FSS}
that there are ``spurious'' modular invariants satisfying
the above constraints but which do not admit an extended
modular S-matrix. But even with this relatively simple
specification we have another problem: Complete classifications
of such modular invariant matrices are known only for
very few models, not much more than $\bbZ_n$ conformal field
theories \cite{D}, $\SUz$ all levels \cite{CIZ2,Kt},
$\SUd$ all levels \cite{G2}, and some classifications
for affine partition functions at low levels \cite{G3}.

Another specification of ``physical'' (but unfortunately
mathematically harder to reach) would be that $Z$ arises from
``the existence of some 2D conformal field theory''.
A promising way of making this precise seems for us to be
the concept of chiral observables as light-cone nets
built in an observable net over 2D Minkowski space \cite{R5}.
As mentioned in Section \ref{t2ext}, Rehren has
shown \cite{R7} that any subfactor $N\subset M$ of our
kind which arises as an extension of a
local factor $N=N(I_\circ)$ of a M\"obius covariant net
$\{N(I)\}$ over $\bbR$ (or equivalently $S^1\setminus\zeta$)
determines an entire 2D conformal field theory over
Minkowski space. The converse direction, however,
is an open problem: Does any 2D conformal field theory
with chiral building blocks containing $\{N(I)\}$
determine a subfactor $N\subset M$ producing the
modular invariant matrix $Z$ which describes the
coupling between left- and right-chiral sectors?
(In particular in the case that the coupling matrix
is \typeii.)
Nevertheless there are partial answers to this question.
First of all the trivial invariants,
$Z_{\la,\mu}=\delta_{\la,\mu}$, are obtained from
the trivial subfactor $N\subset M$ with $M=N$.
Next, any conformal inclusion determines
a subfactor which in turn produces a modular invariant,
being the \typei\ exceptional invariant which arises from
the diagonal invariant of the extended theory,
here the level 1 representation theory of the
larger affine Lie algebra (e.g.\ of $\SOf$ for the
embedding $\SUz_{10}\subset\SOf_1$ as treated above).
The situation is even better for simple current invariants,
which in a sense produce the majority of non-trivial
modular invariants. Simple currents \cite{SY1}
are primary fields
with unit quantum dimension and appear in our framework
as sectors with statistical dimension one, hence its
representatives are automorphisms. They form a closed abelian
group $G$ under fusion which is hence a product of
cyclic groups. Simple currents give rise to modular invariants,
and all such invariants have been classified \cite{GRS,KS}.

If we take generators $[\sigma_i]$ for each cyclic subgroup
$\bbZ_{n_i}$ then we can construct the crossed product subfactor
$N\subset M=N\rtimes G$ whenever we can choose a representative
$\sigma_i$ in each such simple current sector such that we have
exact cyclicity $\sigma_i^{n_i}=\id$ (and not only as sectors).
As we are starting with a chiral quantum field theory
(e.g.\ from loop groups), Rehren's lemma \cite{R3} applies
which states that such a choice is possible if and only if
the statistics phase is an $n_i$-th root of unity, or in the
conformal context if and only if the conformal weight
$h_{\sigma_i}$ is an integer multiple of $1/n_i$.
Sometimes this may only be possible for a simple
current subgroup $H\subset G$, but any non-trivial
subgroup ($H\neq\{0\}$) gives rise to a non-trivial
subfactor and in turn to a modular invariant.
In fact one can check by our methods that all simple
current invariants are realized this way. For example,
for $\SUn_k$ the simple current group is just $\bbZ_n$,
corresponding to weights $k\Lambda_{(j)}$, $j=0,1,...,n-1$.
The conformal dimensions are $h_{k\Lambda_{(j)}}=kj(n-j)/2n$
which allow for extensions except when $n$ is even and $k$
and $j$ are odd. (This reflects the fact that e.g.\
for $\SUz$ there are no D-invariants at odd levels.)
An extension by a simple current subgroup $\bbZ_m\subset\bbZ_n$,
i.e.\ $m$ is a divisor of $n$, is moreover local, if the
generating current (and hence all in the $\bbZ_m$ subgroup)
has integer conformal weight, $h_{k\Lambda_{(q)}}\in\bbZ$,
where $n=mq$. This happens exactly if $kq\in 2m\bbZ$ if $n$
is even, or $kq\in m\bbZ$ if $n$ is odd.
For $\SUz$ this corresponds to the
$\rmD_{\mathrm{even}}$ series whereas the $\rmD_{\mathrm{odd}}$
series are non-local extensions. For $\SUd$, there is a simple
current extension at each level, but only those at $k\in 3\bbZ$
are local. Clearly, the cases with chiral locality match
exactly the \typei\ simple current modular invariants.
Our results imply that the system $\MXMo$ of neutral
morphisms, which is obtained by decomposing $[\a^\pm_\la]$'s
with colour zero mod $m$, carries a non-degenerate braiding.
This nicely reflects a general fact about non-degenerate
extensions of degenerate (sub-) systems conjectured by
Rehren \cite{R1} and proven by M\"uger \cite{Mu}.

For the exceptional modular invariants arising from conformal
inclusions, the corresponding subfactor comes (almost) for free.
A conformal inclusion means that the level 1 representations
of some loop group of a Lie group restrict in a finite manner
to the positive energy representations of a certain embedded
loop group of an embedded (simple) Lie group at some level.
As discussed for the E$_6$ example, a subfactor is obtained
by taking this embedding as a local subfactor in the vacuum
representation. Since the embedding level one theory is always
local, the modular invariant will necessarily be \typei.
For $\SUz$, the modular invariants arising from conformal
embeddings are, besides E$_6$, the E$_8$ and the D$_4$ ones,
corresponding to embeddings $\SUz_{28}\subset(\Gtwo)_1$
and $\SUz_4\subset\SUd_1$, respectively, the latter happens
to be a simple current invariant at the same time.
For $\SUd$, the invariants from conformal embeddings
are $\cD^{(6)}$, $\cE^{(8)}$, $\cE^{(12)}$ and $\cE^{(24)}$,
corresponding to $\SUd_3\subset{\mathit{SO}}(8)_1$,
$\SUd_5\subset{\mathit{SU}}(6)_1$,
$\SUd_9\subset(\rmE_6)_1$, $\SUd_{21}\subset(\rmE_7)_1$,
respectively.

With these techniques we can obtain a huge amount of
modular invariants from subfactors. Nevertheless we still
do not have a systematic procedure to get all physical
invariants. The more problematic cases are typically the
exceptional \typeii\ invariants. We did realize the
E$_7$ invariant of $\SUz$ by some subfactor, namely
we used the existence of a certain
Goodman-de la Harpe-Jones subfactor \cite{GHJ} for this case,
however, this method will not apply to general invariants
of $\SUn$. It seems to follow from Ocneanu's recent announcement
(see his lectures) that there are subfactors realizing all
$\SUd$ modular invariants, but also his methods relying on the
``$\SUd$ wire model'' (as well as on Gannon's classification of
modular invariants) do not solve the general problem.
Nevertheless a large class of exceptional \typeii\ invariants
can be dealt with quite generally, namely those which are
\typeii\ descendants of conformal embeddings. Since
the embedding level 1 theories are typically
(whenever simply laced Lie groups are worked with) $\bbZ_n$
theories, i.e.\ pure simple current theories, the subfactors
producing their modular invariants can be constructed
by simple current methods, and in turn we will obtain the
relevant subfactors for the embedded theories, say $\SUn$.

For a while we will be looking at the so-called $\bbZ_n$
conformal field theories as treated in \cite{D}, which have
$n$ sectors, labelled by $\la=0,1,2,...,n-1$ (mod$\,n$),
obeying $\bbZ_n$ fusion rules, and conformal dimensions
of the form $h_\la=a\la^2/2n$ (mod$\,1$), where $a$ is an integer
mod$\,2n$, $a$ and $n$ coprime and $a$ is even whenever
$n$ is odd. The modular invariants of such models have been
classified \cite{D}. They are labelled by the divisors
$\delta$ of $\tn$, where $\tn=n$ if $n$ is odd and $\tn=n/2$ if
$n$ is even. Explicitly, the modular invariants $Z^{(\delta)}$
are given by
\[
Z^{(\delta)}_{\la,\mu} = \left\{ \begin{array}{ll}
1 & \qquad \mbox{if}\;\;\; \la,\mu=0\,\mbox{mod}\,\alpha
\;\;\;\mbox{and}\;\;\; \mu=\omega(\delta)\la \,\mbox{mod}\,
n/\alpha\,, \\[.4em]
0 & \qquad \mbox{otherwise}\,, \end{array} \right.
\]
where $\a=\mathrm{gcd}(\delta,\tn/\delta)$ so that there
are numbers $r,s\in\bbZ$ such that
$r\tn/\delta\a-s\delta/\a=1$ and then $\omega(\delta)$
is defined as $\omega(\delta)=r\tn/\delta\a+s\delta/\a$.
The trivial invariant corresponds to $\delta=\tn$,
i.e.\ $Z^{(\tn)}=\bfe$ and $\delta=1$ gives the charge
conjugation matrix, $Z^{(1)}=C$.

We now claim that
\begin{equation}\label{Zndiag}
Z^{(\delta)}_{\la,\la} = \left\{ \begin{array}{ll}
1 & \qquad \mbox{if}\;\;\; \la=0\,\mbox{mod}\,\tn/\delta, \\[.4em]
0 & \qquad \mbox{otherwise}\,. \end{array} \right.
\end{equation}
Notice that $\omega(\delta)-1=2s\delta/\alpha$.
Assume first that $\la=x\tn/\delta$, $x\in\bbZ$.
Then clearly $\la=0\,$mod$\,\alpha$ since $\alpha$
divides $\tn/\delta$, and we have
$(\omega(\delta)-1)\la=2sx\tn/\alpha$,
implying $\la=\omega(\delta)\la\,$mod$\,n/\alpha$,
thus $Z^{(\delta)}_{\la,\la}=1$.
Conversely, assume $Z^{(\delta)}_{\la,\la}=1$ so that
$\la=y\alpha$ and $(\omega(\delta)-1)\la=zn/\alpha$
with $y,z\in\bbZ$. This gives
$2sy\delta=zn/\alpha$, hence $2sy=zn/\delta\alpha$.
Now $s$ is coprime to $\tn/\delta\alpha$, and
therefore it follows that $y$ is a multiple of
$\tn/\delta\alpha$ (as we see that $z$ must be
even if $n$ is odd) which implies in fact
$\la=0\,$mod$\,\tn/\delta$.

{}From \erf{Zndiag} we obtain the following trace property
of $Z^{(\delta)}$:
\[
\tr(Z^{(\delta)}) = \epsilon\delta \,, \qquad
\mbox{where}\quad \epsilon= \frac n\tn = \left\{
\begin{array}{cl} 2 & \qquad \mbox{if $n$ is even,}\\[.4em]
1 & \qquad \mbox{if $n$ is odd.} \end{array} \right.
\]
Now suppose that for such a $\bbZ_n$ theory at hand
we have corresponding braided endomorphisms $\rho_\la$
of some \typeiii\ factor $N$, such that their
statistical phases are given by $\E^{2\pi\I h_\la}$
with conformal weights $h_\la$ as above
(as is the case for level 1 loop group theories).
As we are dealing with $\bbZ_n$ fusion rules,
all our morphisms $\rho_\la$ will in fact be automorphisms.
Note that if $n$ is odd then we can always assume that
$\rho_1^n=\id$ as morphisms (and our system can be
chosen as $\{\rho_1^\la\}_{\la=0}^{n-1}$). However,
if $n$ is even, then we cannot choose a representative
of the sector $[\rho_1]$ such that its $n$-th power gives
the identity, nevertheless we can always assume that
$\rho_\epsilon^\tn=\id$.
Thus we have a simple current (sub-) group $\bbZ_\tn$,
for which we can form the crossed product subfactor
$N\subset M=N\rtimes\bbZ_{\tn/\delta}$ for any divisor
$\delta$ of $\tn$. It is quite easy to see that
$N\subset M=N\rtimes\bbZ_{\tn/\delta}$ indeed realizes
$Z^{(\delta)}$: The crossed product by $\bbZ_{\tn/\delta}$
gives the dual canonical endomorphism sector
$[\canr]=[\id]\oplus[\rho_{\epsilon\delta}]\oplus
[\rho_{\epsilon\delta}^2]\oplus\ldots\oplus
[\rho_{\epsilon\delta}^{\tn/\delta-1}]$.
The formula
$\langle\iota\rho_\la,\iota\rho_\mu\rangle=
\langle\canr\rho_\la,\rho_\mu\rangle$ then shows that
the system of $M$-$N$ morphisms is labelled by
$\bbZ_n/\bbZ_{\tn/\delta}\simeq\bbZ_{\epsilon\delta}$,
i.e.\ $\#\MXN=\epsilon\delta$.
Therefore our general theory implies that the modular invariant
arising from $N\subset M=N\rtimes\bbZ_{\tn/\delta}$ has trace
equal to $\epsilon\delta$, and thus must be $Z^{(\delta)}$.
Thus all modular invariants classified in \cite{D} are
realized from subfactors.

It is instructive to apply the above results to descendant
modular invariants of conformal inclusions. Let us consider
the conformal inclusion $\SUf_6\subset{\mathit{SU}}(10)_1$.
The associated modular invariant, which can be found
in \cite{SY1}, reads
\[
Z = \sum_{j\in\bbZ_{10}} |\chi^j|^2
\]
with ${\mathit{SU}}(10)_1$ characters decomposing into
$\SUf_6$ characters as
\[
\begin{array}{ll}
\chi^0 = \chi_{0,0,0} + \chi_{0,6,0} + \chi_{2,0,2} + \chi_{2,2,2}, &
\chi^5 = \chi_{0,0,6} + \chi_{6,0,0} + \chi_{0,2,2} + \chi_{2,2,0},
\\[.4em]
\chi^1 = \chi_{0,0,2} + \chi_{2,4,0} + \chi_{2,1,2}\,, &
\chi^6 = \chi_{4,0,0} + \chi_{0,2,4} + \chi_{1,2,1}\,,
\\[.4em]
\chi^2 = \chi_{0,1,2} + \chi_{2,3,0} + \chi_{3,0,3}\,, &
\chi^7 = \chi_{3,0,1} + \chi_{1,2,3} + \chi_{0,3,0}\,,
\\[.4em]
\chi^3 = \chi_{1,0,3} + \chi_{3,2,1} + \chi_{0,3,0}\,, &
\chi^8 = \chi_{0,3,2} + \chi_{2,1,0} + \chi_{3,0,3}\,,
\\[.4em]
\chi^4 = \chi_{0,0,4} + \chi_{4,2,0} + \chi_{1,2,1}\,, &
\chi^9 = \chi_{2,0,0} + \chi_{0,4,2} + \chi_{2,1,2}\,.
\end{array}
\]
As usual, this invariant can be realized from the conformal
inclusion subfactor
\[
N = \pi^0({\mathit{L}}_I{\mathit{SU}}(4))'' \subset
\pi^0({\mathit{L}}_I{\mathit{SU}}(10))'' = M_+ \,,
\]
with $\pi^0$ denoting the level 1 vacuum representation
of $\mathit{LSU}(10)$. The dual canonical endomorphism
sector corresponds to the vacuum block,
\[
[\canr_+] = [\la_{0,0,0}]\oplus[\la_{0,6,0}]\oplus
[\la_{2,0,2}]\oplus[\la_{2,2,2}]\,.
\]
Proceeding with $\a$-induction
$\la_{p,q,r}\mapsto\a^\pm_{+;p,q,r}\in\End(M_+)$,
it is a straightforward calculation that the graphs
describing left multiplication by fundamental generators
$[\a^\pm_{+;1,0,0}]$ and $[\a^\pm_{+;0,1,0}]$ (which is
the same as right multiplication by $[\la_{1,0,0}]$ and
$[\la_{0,1,0}]$, respectively) on the system of $M_+$-$N$
sectors gives precisely the graphs found by Petkova and
Zuber \cite[Figs.\ 1 and 2]{PZp} by their more empirical
procedure to obtain graphs with spectrum matching the diagonal
part of some given modular invariant. In our framework,
the graph \cite[Fig.\ 1]{PZp} obtains the following meaning:
Take the outer wreath, pick a vertex with 4-ality 0 and
label it by $[\iota_+]\equiv[\tau_0\iota_+]$,
where $\iota_+:N\hookrightarrow M_+$ denotes the
injection homomorphism, as usual.
Going around in a counter-clockwise direction the vertices
will then be the marked vertices labelled by the $\bbZ_{10}$
sectors $[\tau_1\iota_+]$, $[\tau_2\iota_+]$, .... ,
$[\tau_9\iota_+]$ of $\mathit{SU}(10)_1$.
Passing to the next inner wreath
the 4-ality 1 vertex adjacent to $[\iota_+]$ is then
the sector $[\a^\pm_{+;1,0,0}\iota_+]=[\iota_+\la_{1,0,0}]$,
and the others its $\bbZ_{10}$ translates.
Similarly the inner wreath consists
of the $\bbZ_{10}$ translates of $[\iota_+\la_{0,1,0}]$.
The remaining two vertices in the center correspond to
subsectors of the reducible $[\iota\la_{1,1,0}]$ and
$[\iota\la_{0,1,1}]$. The graph itself then represents
left (right) multiplication by $[\a^\pm_{+;1,0,0}]$
($[\la_{1,0,0}]$).

As for $\mathit{LSU}(10)$ at level 1 we are in fact dealing
with a $\bbZ_n$ conformal field theory, we have $n=10$ and
$\tn=5$, we thus know that there are only two modular
invariants: The diagonal one which in restriction to
$\mathit{LSU}(4)$ gives exactly the above \typei\ invariant
$Z\equiv Z^{(5)}$, but there is also the charge conjugation
invariant $Z^{(1)}$, written as
\[
Z^{(1)} = \sum_{j\in\bbZ_{10}} \chi^j (\chi^{-j})^* \,.
\]
Whereas $Z^{(5)}$ can be thought of as the trivial
extension $M_+\subset M_+$, the conjugation invariant $Z^{(1)}$
can be realized from the crossed product
$M_+\subset M=M_+\rtimes\bbZ_5$ which has dual canonical
endomorphism sector
\[
[\canr^\ext]=[\tau_0]\oplus[\tau_2]\oplus[\tau_4]
\oplus[\tau_6]\oplus[\tau_8] \,.
\]
So far we have considered the situation on the ``extended level'',
but we may now descend to the level of $\SUf_6$ sectors and
characters. Namely we may consider the subfactor
$N\subset M=M_+\rtimes\bbZ_5$. Its dual canonical endomorphism
sector $[\canr]$ is obtained by $\sigma$-restriction of $[\canr^\ext]$
which can now be read off from the character decomposition,
\[
\begin{array}{rl}
[\canr]&=\bigoplus_{j=0}^4 [\sigma_{\tau_{2j}}]
=[\la_{0,0,0}]\oplus[\la_{0,6,0}]\oplus[\la_{2,0,2}]
\oplus[\la_{2,2,2}]\oplus [\la_{0,1,2}] \\[.4em]
& \quad
\oplus[\la_{2,3,0}]\oplus[\la_{3,0,3}]\oplus
[\la_{0,0,4}]\oplus[\la_{4,2,0}]\oplus[\la_{1,2,1}]\oplus
[\la_{4,0,0}]\oplus[\la_{0,2,4}] \\[.4em]
& \quad
\oplus[\la_{1,2,1}]\oplus
[\la_{0,3,2}]\oplus[\la_{2,1,0}]\oplus[\la_{3,0,3}]
\,. \end{array}
\]
This subfactor produces the conjugation invariant
$Z^{(1)}$ written in $\SUf_6$ characters which is
the same as taking the original $\SUf_6$ conformal
inclusion invariant and conjugating on the level of
the $\SUf_6$ characters. Note that this invariant has
only 16 diagonal entries.

Also note that we will still have entries
$Z_{\la,\mu}\ge 2$, for instance the diagonal
entry corresponding to the weight $(2,1,2)$ is $2$
as $|\chi_{2,1,2}|^2$ appears in $\chi^1(\chi^9)^*$
and in $\chi^9(\chi^1)^*$. Hence the system of
$M$-$M$ sectors will have non-commutative fusion rules
(as had the $M_+$-$M_+$ system).
When passing from $M_+$ to $M=M_+\rtimes\bbZ_5$,
the $M_+$-$N$ system will change
to the $M$-$N$ system in such a way
that all sectors which are translates
by $\tau_{2j}$, $j=0,1,2,3,4$, have to be identified,
and similarly fixed points split. Thus our new system of
$M$-$N$ morphisms will be some kind of orbifold of the old one.
To see this, we first recall that all the irreducible
$M_+$-$N$ morphisms are of the form $\beta\iota_+$ with
$\beta\in\MXMppm$. To such an irreducible $M_+$-$N$ morphism
$\beta\iota_+$ we can now associate an $M$-$N$ morphism
$\iota^\ext\beta\iota_+$ which may no longer be irreducible;
here $\iota^\ext$ is the injection homomorphism
$M_+\hookrightarrow M$. Then the reducibility can be
controlled by Frobenius reciprocity as we have
\[
\langle\iota^\ext\beta\iota_+,\iota^\ext\beta'\iota_+\rangle=
\langle\canr^\ext\beta\iota_+,\beta'\iota_+\rangle \,,
\]
and $\canr^\ext=\co\iota ^\ext\iota^\ext$.
Carrying out the entire computation we find that there
are 16 $M$-$N$ sectors, and the
right multiplication by $[\la_{1,0,0}]$ is displayed
graphically as in \fig{orbPZ1}.
\begin{figure}[htb]
\begin{center}
\unitlength 0.8mm
\begin{picture}(80,110)
\thinlines
\multiput(30,50)(40,0){2}{\circle*{4}}
\multiput(10,50)(40,0){2}{\circle{4}}
\multiput(10,50)(40,0){2}{\circle{3}}
\multiput(10,50)(40,0){2}{\circle*{2}}
\put(40,70){\circle{4}}
\put(40,70){\circle*{3}}
\put(40,30){\circle{4}}
\put(40,30){\circle{3}}
\put(40,30){\circle{2}}
\put(40,30){\circle*{1}}
\multiput(40,80)(0,5){5}{\circle{4}}
\multiput(40,80)(0,5){5}{\circle{3}}
\multiput(40,80)(0,5){5}{\circle{2}}
\multiput(40,80)(0,5){5}{\circle*{1}}
\multiput(40,0)(0,5){5}{\circle{4}}
\multiput(40,0)(0,5){5}{\circle*{3}}
\thicklines
\path(30,50)(40,70)(50,50)(40,30)(30,50)
\path(10,50)(40,80)(70,50)(40,20)(10,50)
(40,85)(70,50)(40,15)(10,50)
(40,90)(70,50)(40,10)(10,50)
(40,95)(70,50)(40,5)(10,50)
(40,100)(70,50)(40,0)(10,50)
\path(9.5,50.5)(39.5,70.5)
\path(10.5,49.5)(40.5,69.5)
\path(9.5,49.5)(39.5,29.5)
\path(10.5,50.5)(40.5,30.5)
\path(39.5,30.5)(69.5,50.5)
\path(40.5,29.5)(70.5,49.5)
\path(39.5,69.5)(69.5,49.5)
\path(40.5,70.5)(70.5,50.5)
\end{picture}
\end{center}
\caption{Graph $G_1$ associated to the conjugation invariant
of the conformal inclusion $\SUf_6\subset{\mathit{SU}}(10)_1$}
\label{orbPZ1}
\end{figure}
Here the 4-alities 0,1,2,3 of the vertices are indicated
by solid circles of decreasing size. The $[\iota]$
vertex (with $\iota=\iota^\ext\iota_+$ denoting the
injection homomorphism $N\hookrightarrow M$ of the
total subfactor $N\subset M=M_+\rtimes \bbZ_5$)
is the 4-ality 0 vertex in the center of the picture,
and the 4-ality 1 vertex above corresponds to
$[\iota\la_{1,0,0}]$. Each group of five vertices on
the top and the bottom of the picture arise from the
splitting of the two central vertices of the graphs
in \cite{PZp} as they are $\bbZ_5$ fixed points.
That our orbifold graph inherits the 4-ality of the
original graph is due to the fact that all entries in
$[\canr]$ have 4-ality zero which in turn comes from the
fact that all even marked vertices (corresponding to
the subgroup $\bbZ_5\subset\bbZ_{10}$) of the graph of
Petkova and Zuber have 4-ality zero.
We also display the graph corresponding to the second
fundamental representation, namely the right multiplication
by $[\la_{0,1,0}]$ in \fig{orbPZ1}.
\begin{figure}[htb]
\begin{center}
\unitlength 0.8mm
\begin{picture}(120,110)
\thinlines
\multiput(60,50)(40,0){2}{\circle*{4}}
\multiput(40,50)(40,0){2}{\circle{4}}
\multiput(40,50)(40,0){2}{\circle{3}}
\multiput(40,50)(40,0){2}{\circle*{2}}
\put(70,70){\circle{4}}
\put(70,70){\circle*{3}}
\put(70,30){\circle{4}}
\put(70,30){\circle{3}}
\put(70,30){\circle{2}}
\put(70,30){\circle*{1}}
\multiput(70,80)(0,5){5}{\circle{4}}
\multiput(70,80)(0,5){5}{\circle{3}}
\multiput(70,80)(0,5){5}{\circle{2}}
\multiput(70,80)(0,5){5}{\circle*{1}}
\multiput(70,0)(0,5){5}{\circle{4}}
\multiput(70,0)(0,5){5}{\circle*{3}}
\thicklines
\path(40,50)(60,50)
\path(80,50)(100,50)
\put(70,21){\arc{84.853}{3.927}{5.498}}
\put(70,19){\arc{84.853}{3.927}{5.498}}
\put(70,81){\arc{84.853}{0.785}{2.356}}
\put(70,79){\arc{84.853}{0.785}{2.356}}
\put(65,25){\arc{14.142}{5.498}{0.785}}
\put(62.5,22.5){\arc{21.213}{5.498}{0.785}}
\put(60,20){\arc{28.284}{5.498}{0.785}}
\put(57.5,17.5){\arc{35.355}{5.498}{0.785}}
\put(55,15){\arc{42.426}{5.498}{0.785}}
\put(65,75){\arc{14.142}{5.498}{0.785}}
\put(62.5,77.5){\arc{21.213}{5.498}{0.785}}
\put(60,80){\arc{28.284}{5.498}{0.785}}
\put(57.5,82.5){\arc{35.355}{5.498}{0.785}}
\put(55,85){\arc{42.426}{5.498}{0.785}}
\put(60,50){\arc{61}{1.571}{4.712}}
\put(60,50){\arc{59}{1.571}{4.712}}
\put(60,50){\arc{71}{1.571}{4.712}}
\put(60,50){\arc{69}{1.571}{4.712}}
\put(60,50){\arc{81}{1.571}{4.712}}
\put(60,50){\arc{79}{1.571}{4.712}}
\put(60,50){\arc{91}{1.571}{4.712}}
\put(60,50){\arc{89}{1.571}{4.712}}
\put(60,50){\arc{101}{1.571}{4.712}}
\put(60,50){\arc{99}{1.571}{4.712}}
\path(60,0.5)(70,0.5)
\path(60,-0.5)(70,-0.5)
\path(60,5.5)(70,5.5)
\path(60,4.5)(70,4.5)
\path(60,10.5)(70,10.5)
\path(60,9.5)(70,9.5)
\path(60,15.5)(70,15.5)
\path(60,14.5)(70,14.5)
\path(60,20.5)(70,20.5)
\path(60,19.5)(70,19.5)
\path(60,80.5)(70,80.5)
\path(60,79.5)(70,79.5)
\path(60,85.5)(70,85.5)
\path(60,84.5)(70,84.5)
\path(60,90.5)(70,90.5)
\path(60,89.5)(70,89.5)
\path(60,95.5)(70,95.5)
\path(60,94.5)(70,94.5)
\path(60,100.5)(70,100.5)
\path(60,99.5)(70,99.5)
\put(90,50){\arc{41}{4.712}{1.571}}
\put(90,50){\arc{39}{4.712}{1.571}}
\path(70,30.5)(90,30.5)
\path(70,29.5)(90,29.5)
\path(70,70.5)(90,70.5)
\path(70,69.5)(90,69.5)
\end{picture}
\end{center}
\caption{Graph $G_2$ associated to the conjugation invariant
of the conformal inclusion $\SUf_6\subset{\mathit{SU}}(10)_1$}
\label{orbPZ2}
\end{figure}
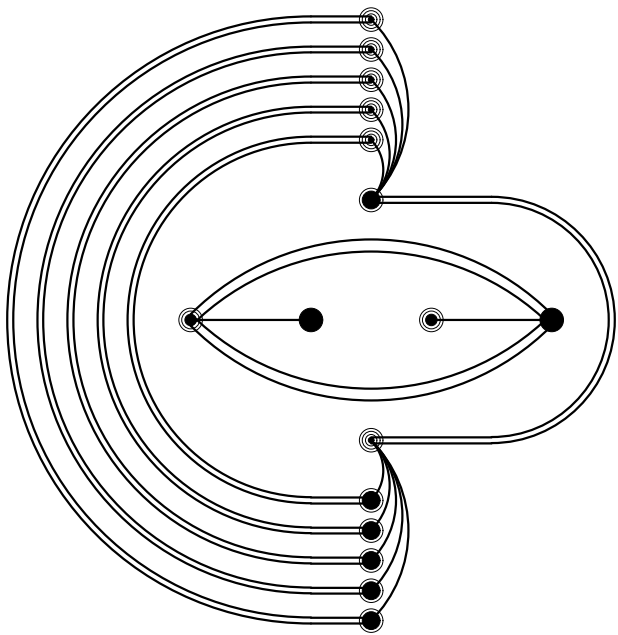

The conformal inclusion $\SUd_5\subset{\mathit{SU}}(6)_1$
can be treated along the same lines. The associated $\SUd_5$
modular invariant, i.e.\ the one which is the specialization
of the diagonal ${\mathit{SU}}(6)_1$ invariant,
\[
Z=\sum_{j\in\bbZ_6} |\chi^j|^2 \,,
\]
with ${\mathit{SU}}(6)_1$ characters decomposing
in $\SUd_5$ variables as
\[
\begin{array}{lll}
  \chi^0 = \chi_{0,0} + \chi_{2,2} \,,\quad
& \chi^1 = \chi_{2,0} + \chi_{2,3} \,,\quad
& \chi^2 = \chi_{2,1} + \chi_{0,5} \,, \\[.4em]
  \chi^3 = \chi_{3,0} + \chi_{0,3} \,,
& \chi^4 = \chi_{1,2} + \chi_{5,0} \,,
& \chi^5 = \chi_{0,2} + \chi_{3,2} \,,
\end{array}
\]
is labelled by the graph $\cE^{(8)}$. Besides this diagonal
invariant $Z\equiv Z^{(3)}$, the extended ${\mathit{SU}}(6)_1$
theory, being a $\bbZ_6$ theory, possesses only the conjugation
invariant $Z^{(1)}=\sum_{j\in\bbZ_6} \chi^j(\chi^{-j})^*$,
corresponding to the divisors 3 and 1 of 3,
respectively. Writing again the conformal inclusion subfactor
as $N\subset M_+$, the conjugation invariant can be
realized from the extension $N\subset M=M_+\rtimes \bbZ_3$
with canonical endomorphism sector
\[
[\canr] = [\la_{0,0}] \oplus [\la_{2,2}] \oplus [\la_{2,1}]
\oplus [\la_{0,5}] \oplus [\la_{1,2}] \oplus [\la_{5,0}] \,,
\]
which arises as $\canr=\sigma_{\canr^\ext}$ where
$[\canr^\ext]=[\tau_0]\oplus[\tau_2]\oplus[\tau_4]$.
Whereas the $M_+$-$N$ system is labelled by the vertices of
the graph $\cE^{(8)}$ and can be given by $\{\beta\iota\}$
where $\beta$ runs through the chiral $M_+$-$M_+$ system
determined in \cite[Subsect.\ 2.3 (iv)]{BE2}, the
$M$-$N$ system will now be obtained from this one by
identification of all $\bbZ_3$ translations
(corresponding to the vertices labelled by
$[\a_{(0,0)}]$, $[\a_{(5,5)}]$ and $[\a_{(5,0)}]$
in \cite[Fig.\ 11]{BE2}). We have no fixed points here
so that the 12 vertices of $\cE^{(8)}$ collapse to
4 vertices, and it is easy to see that the new
$M$-$N$ fusion graph is exactly the graph ${\cE^{(8)}}^*$
in the list of Di Francesco and Zuber
(see Zuber's lectures or \cite{BPPZ}).
Note that this time the orbifold graph (${\cE^{(8)}}^*$)
looses the triality of the original graph ($\cE^{(8)}$)
because the even marked vertices (corresponding to the
subgroup $\bbZ_3\subset\bbZ_6$) of $\cE^{(8)}$
are not exclusively of colour zero.

This way we understand why the descendants of modular
invariants of conformal inclusions (where the extended
theory has $\bbZ_n$ fusion rules) are in fact labelled
by orbifold graphs of the graph labelling the original,
block-diagonal conformal inclusion invariant, and
why the conjugation invariant corresponds to the
maximal $\bbZ_\tn$ orbifold.

In the above examples, the trivial and conjugation invariant
of the extended theory still remained distinct when written
in terms of the $\SUf_6$ characters. This need not be the
case in general. Let us look at a familiar modular invariant
of $\SUd$ at level 9, namely
\[
Z_{\cE^{(12)}} = | \chi_{0,0} + \chi_{9,0} +
\chi_{0,9} + \chi_{4,1} + \chi_{1,4} + \chi_{ 4,4}|^2
+ 2\, |\chi_{2,2} + \chi_{5,2} + \chi_{2,5} |^2 \,,
\]
which arises from the conformal embedding
$\SUd_9\subset(\rmE_6)_1$. Now $\rmE_6$
at level 1 gives a $\bbZ_3$ theory and in terms of
the extended characters the above invariant is the
trivial extended invariant
\[
Z_{\cE_1^{(12)}} = |\chi^0|^2+|\chi^1|^2+|\chi^2|^2\,,
\]
using obvious notation. Here both the $(\rmE_6)_1$
characters $\chi^1$ and $\chi^2$ specialize to
$\chi_{2,2} + \chi_{5,2} + \chi_{2,5}$ in terms of
$\SUd_9$ variables. Let $N\subset M_+$ denote the conformal
inclusion subfactor obtained by analogous means as in the
previous example. It has been treated in \cite{BE3}
and produces the graph $\cE^{(12)}_1$ of the list of
Di Francesco and Zuber as chiral fusion graphs ---
and in turn as $M_+$-$N$ fusion graph,
thanks to chiral locality.

Corresponding to the two divisors
3 and 1 of 3, we know that besides the trivial there
is only the conjugation invariant of our $\bbZ_3$ theory.
It is given as
\[
Z_{\cE_2^{(12)}} = |\chi^0|^2+\chi^1(\chi^2)^*
+\chi^2(\chi^1)^*
\]
but this distinct invariant restricts to the same
invariant $Z_{\cE^{(12)}}$ when specialized to
$\SUd_9$ variables. Nevertheless we will obtain a different
subfactor $N\subset M$ since the conjugation invariant
of our $\bbZ_3$ theory is realized from the extension
$M_+\subset M=M_+\rtimes\bbZ_3$. In particular,
the subfactor $N\subset M$ has dual canonical
endomorphism sector
\[
[\canr] = [\la_{0,0}] \oplus [\la_{9,0}] \oplus
[\la_{0,9}] \oplus [\la_{4,1}] \oplus [\la_{1,4}] \oplus
[\la_{ 4,4}] \oplus 2 [\la_{2,2}] \oplus 2 [\la_{5,2}]
\oplus 2[\la_{2,5}] \,,
\]
determined by $\sigma$-restriction of
\[
[\canr^\ext] = [\tau_0] \oplus [\tau_1] \oplus [\tau_2 ] \,.
\]
As before, the $M$-$N$ system can be obtained from the $M_+$-$N$
system by dividing out the cyclic symmetry carried by
$[\canr^\ext]$. In terms of graphs, the cyclic $\bbZ_3$ symmetry
corresponds to the three wings of the graph $\cE^{(12)}_1$ which
are transformed into each other by translation through the
$[\tau_j]$'s, and dividing out this symmetry gives exactly the
graph $\cE^{(12)}_2$ as the wings are identified whereas each
vertex on the middle axis splits into three nodes of identical
Perron-Frobenius weight. This way we understand the graph
$\cE^{(12)}_2$ as the label for the conjugation invariant
$Z_{\cE_2^{(12)}}$ of $Z_{\cE_1^{(12)}}$ which accidentally
happens to be the same as the selfconjugate $Z_{\cE^{(12)}}$ when
specialized to $\SUd_9$ variables.

Though here the same modular invariant, possessing a second
interpretation as its own conjugation, gave rise to two different
graphs, it often happens that an exceptional self-conjugate
invariant is labelled by only one and the same graph
which is its own orbifold. The very simplest case is the
conformal inclusion $\SUz_4\subset\SUd_1$, giving rise
to the D$_4$ invariant which is self-conjugate for $\SUz$
though the non-specialized diagonal $\SUd_1$ invariant is not.
We could proceed as above, passing from the conformal inclusion
subfactor $N\subset M_+$ to $N\subset M=M_+\rtimes\bbZ_3$,
collapsing the $M_+$-$N$ fusion graph D$_4$ into its
$\bbZ_3$ orbifold. However, identifying the three external
vertices and splitting the $\bbZ_3$ fixed point into 3
nodes gives us again D$_4$: The Dynkin diagram D$_4$ is
its own $\bbZ_3$ orbifold.

\vspace{0.5cm}\addtolength{\baselineskip}{-2pt}
\begin{footnotesize}
\noindent{\it Acknowledgement.}
We would like to thank J.\ Fuchs, T.\ Gannon, C.\ Schweigert
and J.-B.\ Zuber for helpful comments on an earlier
version of the manuscript.
\end{footnotesize}
\vspace{0.5cm}



\newcommand\bitem[2]{\bibitem{#1}{#2}}

\def\aam              {Acta Appl.\ Math. }
\def\aip              {Ann.\ Inst.\ H.\ Poincar\'e (Phys.\ Th\'eor.) }
\def\cmp              {Com\-mun.\ Math.\ Phys. }
\def\duke             {Duke Math.\ J. }
\def\ijm              {Intern.\ J. Math. }
\def\jfa              {J.\ Funct.\ Anal. }
\def\jmp              {J.\ Math.\ Phys. }
\def\lmp              {Lett.\ Math.\ Phys. }
\def\rmp              {Rev.\ Math.\ Phys. }
\def\inv              {Invent.\ Math. }
\def\mpl              {Mod.\ Phys.\ Lett. }
\def\nup              {Nucl.\ Phys. }
\def\nupp             {Nucl.\ Phys.\ (Proc.\ Suppl.) }
\def\adma             {Adv.\ Math. }
\def\physa            {Physica \textbf{A} }
\def\ijmp             {Int.\ J.\ Mod.\ Phys. }
\def\jp               {J.\ Phys. }
\def\fdp              {Fortschr.\ Phys. }
\def\pl               {Phys.\ Lett.}
\def\rims             {Publ.\ RIMS, Kyoto Univ. }


\bibliographystyle{amsalpha}

\end{document}